\documentclass[referee, sn-aps]{sn-jnl}
\usepackage[T1]{fontenc}
\usepackage{amsmath}
\usepackage{amssymb}
\usepackage{amsthm}
\usepackage{graphicx}
\usepackage{tensor}
\usepackage{mathrsfs}

\theoremstyle{thmstyletwo}
\newtheorem{proposition}{Proposition}

\theoremstyle{remark}

\newtheorem{rem}{Remark}

\begin{document}
\title{ Dynamics of spinning test bodies in the Schwarzschild space-time: reduction and
circular orbits}
\author{Ivan Bizyaev }
\affil{Ural Mathematical Center, Udmurt State University, ul. Universitetskaya 1, 426034 Izhevsk, Russia}
\email{ bizyaevtheory@gmail.com}

\abstract{
This paper investigates the motion of a rotating test body
in the Schwarzschild space-time.
		Previously, it was shown that this problem reduces to investigating a two-dimensional
Poincar\'{e} map. The paper presents a detailed analysis of
bifurcations of periodic solutions
using this map.
		In the Poincar\'{e} map, as the energy of the
body increases, one can observe two pitchfork bifurcations that follow one after the other:
a supercritical and a subcritical one. This gives rise to five fixed points in the Poincar\'{e} map.
		 In addition, new circular orbits are found for which the total angular momentum
is not parallel to the angular momentum of the test body.
		 For these circular orbits, the radial coordinate satisfies the condition
$r \geqslant 3$ (in units of the mass of a black hole).
		 For values of the total angular momentum of the test body that corresponds to neutron stars
or black holes, these  asymmetric circular orbits turn out to be unstable.}

	
		\maketitle
	\tableofcontents

	\section{Introduction}
	In general relativity theory, the angular momentum of a test body can have
a considerable influence on its trajectory.
	The method of multipole expansion allows one to obtain a closed finite-dimensional system
for the motion of this body.
	In the first order (pole), the trajectory of the body is a geodesic.
In the next order (pole-dipole),
	the motion of the body is described by the Mathisson\,---\,Papapetrou\,---\,Dixon (MPD) equations
\cite{Mathisson1937, Papapetrou1951, Dixon1964}.
	These equations govern the evolution of the linear momentum $p^\alpha$ and the angular
momentum $S^{\alpha \beta}$ of the body,
	but not the evolution of velocity $u^\alpha$. 	In order to eliminate this
ambiguity, we will use the Tulczyjew condition \cite{Tulczyjew1959}
 and the parameter along the trajectory from Ref. \cite{Ehlers1977} (see Section \ref{sec_tul}).

	The MPD equations reduce to investigating a Hamiltonian system
with six degrees of freedom \cite{Witzany2019}.
For the Schwarzschild space-time, in the standard coordinates
$(t, r, \theta, \varphi)$,
this system has an energy integral
	$\mathcal{E}$ and three additional integrals which define an analog of the ``Euclidean''  total angular momentum vector $\vec{Q}=(Q_1,Q_2, Q_3)$.
	The Poisson bracket between the components of the vector $\vec{Q}$ forms the algebra
$so(3)$. Therefore, if $\vec{Q} \not= 0$, a reduction by three degrees of freedom
is possible since three integrals $\mathcal{E}$, $Q_3$ and $(\vec{Q}, \vec{Q})$ are in involution with each other.
	Given the Tulczyjew condition, this problem reduces to analysis of a two-dimensional
Poincar\'{e} map.

As shown in this paper, if the  total angular momentum vector is zero ($\vec{Q} = 0 $),
then the MPD equations admit a reduction by four degrees of freedom.
In this case, the equations of motion on the invariant manifold given by the Tulczyjew
condition are integrable by quadratures (see Section \ref{sectionQ_zero}). 	

The Poincar\'{e} map that arises in the case $\vec{Q} \not= 0$ was
investigated previously in Ref. \cite{Hartl2003, Suzuki1997}.
There it was shown that in the general case the MPD equations are
nonintegrable by quadratures since in the Poincar\'{e} map
one can clearly see chaotic trajectories for which the leading Lyapunov exponent has
positive values. In Ref. \cite{Zelenka2020}, it is shown that resonances and chaos
can be found for astrophysical values of the angular momentum of the test body.
Recently, there has been growing interest in analyzing the motion of a
rotating test body. This is due to the fact that this problem is central to the modeling of
extreme mass-ratio inspirals (EMRIs) for space-based gravitational-wave observatories
(for details, see \cite{Zelenka2020}) such as LISA, TianQin and Taiji.

In this paper,  a reduction by the symmetry group $SO(3)$ is
explicitly performed for the MPD equations in the Schwarzschild space-time.
The reduced variables used are related to the relativistic Andoyer variables introduced in
Ref. \cite{Ramond2025} (see Section \ref{eqCanonic}).
Note that the equations obtained in Ref. \cite{Ramond2025} in reduced variables are not the
original nonlinear equations, but the nonlinear equations linearized with respect to the angular momentum.
 Another analog of the Andoyer variables was introduced in Ref. \cite{Witzany2019}
 (see Section \ref{eqCanonic}).
The variables used in Ref. \cite{Witzany2019} do not define reduction by the symmetry group
$SO(3)$. Hence, the Hamiltonian obtained in the above-mentioned variables
depends on a larger number of variables (angle $\theta$ and linear momentum $P_\theta$) compared to
the Hamiltonian obtained in this paper after reduction (see Section \ref{sectionRedQneq0}). 

In this paper we numerically construct a Poincar\'{e} map.
It agrees with the map from Refs. \cite{Suzuki1997, Zelenka2020}, which was constructed after
numerical integration of the original MPD equations. It is shown that,
as the energy of the body increases, two pitchfork bifurcations (a supercritical and a
subcritical one) following one after the other can be observed in the
Poincar\'{e} map. This gives rise to five fixed points or five periodic
solutions of the reduced system in the Poincar\'{e} map. In addition, we construct
 stable and unstable branches of the asymptotic manifolds for the unstable fixed points of the Poincar\'{e} map which intersect each other
transversely. To numerically construct the asymptotic manifolds, we use the algorithm suggested in
Ref. \cite{LiUnstable}.

 The fixed points of the reduced system describe a motion where
 the trajectory of the body in space is a circle which has the radial coordinate
 $r={\rm const}$ and the angle $\theta={\rm const}$ (generally speaking,
 $\theta \neq \dfrac{\pi}{2}$ is possible). Throughout this paper, such trajectories will
 be called {\it circular orbits}.
Previously, circular orbits in which the total angular momentum is parallel to the angular
momentum of the body were analyzed in Refs. \cite{Suzuki1998, Jefremov2015, Hackmann2014}.
In this paper we find new circular orbits in which
the total angular momentum is not parallel to the angular momentum of the body.
We will call the orbits from the first case {\it symmetric} (they lie on the line
of fixed symmetry points), and those from the second case, {\it asymmetric}.
For values of the total angular momentum of the test body that corresponds to neutron stars
or black holes, these asymmetric circular orbits turn out to be unstable.

For the geodesics of a Schwarzschild metric that are circular orbits, the
dimensionless radial coordinate expressed in units of the mass of a black hole changes \cite{chandrasekhar1998mathematical} in the interval $r\in(3,+\infty)$.
For the MPD equations, the spin-orbital interaction causes the circular orbits
to move closer to the event horizon \cite{Suzuki1997, Jefremov2015}.
However, for the MPD equations the velocity $u^\alpha$ can fall outside the light cone
(i.e., it can become space-like or isotropic). In this case, the pole-dipole approximation
is obviously violated, and the trajectory becomes nonphysical \cite{Semerak1, HojmanHojman}.
For symmetric circular orbits, this is observed
as the circular orbits approach the event horizon in region $r<3$.
For asymmetric circular orbits, it turns out that $r \geqslant 3 $ and the velocity $u^\alpha$
is a time-like vector (see Section \ref{nonsymmetricOrbit}).

For a rotating test body with an energy-momentum tensor satisfying some natural assumptions, it has been
proved in \cite{Schattner1979a, Schattner1979b} that the Tulczyjew condition is
satisfied by the unique point inside the body. Therefore, this paper
considers the MPD equations only for the Tulczyjew condition. In the literature, use is also made
of other conditions (see, e.g., \cite{Costa2015, Witzany2019}).
 The reduction of the MPD equations for the Mathisson–Pirani and Okahashi–Kyrian–Semer\'ak
 conditions to Hamiltonian form is described in detail in Ref. \cite{Witzany2019}.
 Ref. \cite{Barausse2009} uses the Newton-Wigner condition, for which a Hamiltonian is obtained to
 linear order in the angular momentum $S^{\alpha \beta}$.

Of special note is the fact that, in the case of a rotating test body for which the
Mathisson–Pirani condition is used in the Schwarzschild metric, nonequatorial circular orbits
(i.e., $\theta \neq \dfrac{\pi}{2}$) are described in Refs. \cite{Plyatsko1982, Plyatsko2012}.
Comparison of these orbits with the asymmetric circular
orbits obtained in this paper is a problem in its own right, and so it is not addressed here.
Circular orbits in the equatorial plane are compared in Ref. \cite{Timogiannis2021}
for different additional conditions.

\section{The notation and variables used}

	In this paper we choose a geometric system of units in which the gravitational constant
and the velocity of light are equal to unity and the following notation is used.
	 The Latin indices of tensors range from 1 to 3, and the Greek indices
of tensors range from 0 to 3. Summation over repeating indices is implied. The signature of the metric is $(-,+,+,+)$.
	We denote the antisymmetric Levi-Civita symbol by $\varepsilon_{\alpha \beta \mu \nu}$
for which we set $\varepsilon_{0123}=1$. The covariant derivative is denoted by
$\nabla_\alpha$ and the curvature tensor is defined by $[\nabla_\alpha, \nabla_\beta]A^\mu=\tensor{R}{^\mu_\nu_\alpha_\beta}A^\nu$.

 Three-dimensional vectors will be denoted by a short arrow, i.e.,
$\vec{X}=(X_1, X_2, X_3)$, $\vec{Y}=(Y_1, Y_2, Y_3)$, etc. The scalar and vector products
of the Euclidean vectors
will be denoted by $(\vec{X}, \vec{Y})$ and $\vec{X} \times \vec{Y}$, respectively. The
symbol $\hat{X}_{ij}$ will denote the components of
the following matrix constructed from the Euclidean vector $\vec{X}$:
 $$
\hat{X}_{ij} = \left(
\begin{array}{ccc}
0 & X_3 & -X_2 \\
-X_3 & 0 & X_1 \\
X_2 & -X_1 & 0
\end{array}
\right).
 $$

 In this paper, different sets of variables will be used. The relation between different variables can be illustrated by the following scheme:
	 $$
	 \overset{{\bf I}}{( S^{\alpha \beta}, p_\alpha, y^\alpha )} \xrightarrow{\eqref{eq_Sp}, \eqref{eq_SAB}} \overset{{\bf II}}{(\vec{L}, \vec{M}, P_\alpha, y^\alpha)} \xrightarrow{\eqref{eq_EZ0}} \overset{{\bf III}}{(\vec{E}, \vec{Z}, P_r, r)},
	 $$
	 where the numbers above the arrows between variables denote the numbers of the formulae
in the text which define a transition between them, and the Roman numerals denote the
numbers of the sets of variables.

	 Each set of variables has its advantages.
	 The set of variables ${\bf I}$ is convenient in that the Tulczyjew condition and
the Hamiltonian take the simplest form. For the set of variables ${\bf II}$, the Poisson
bracket has the simplest form. Both of these sets of variables have been used
more than once, for example, in Refs. \cite{Tauber1988, Witzany2019, Ramond2025}.
	 The set of variables ${\bf III}$ is new. Its advantage is that, for the
Schwarzschild metric, it defines the set of reduced variables. The Poisson bracket between
the components of the vectors $\vec{E}$ and $\vec{Z}$ forms the algebra $so(3)$ and $so(2,1)$,
respectively.
	
	\section{The Mathisson\,---\,Papapetrou\,---\,Dixon equations}
	\label{sectionMP}
	In general relativity theory, space-time is considered as a
four-dimensional manifold $\mathcal{M}^4$. On this manifold, we denote the local
coordinates as $(x^{\alpha})=(x^0,\dots, x^3)$ and represent the interval as
	 $$
	 ds^2=g_{\alpha \beta}dx^{\alpha} dx^{\beta},
	 $$
	 where $\boldsymbol{g}=|| g_{\alpha \beta} ||$ is a metric that satisfies the
Einstein equations.
	
\begin{figure}[!ht]
\begin{center}
\includegraphics[scale=0.8]{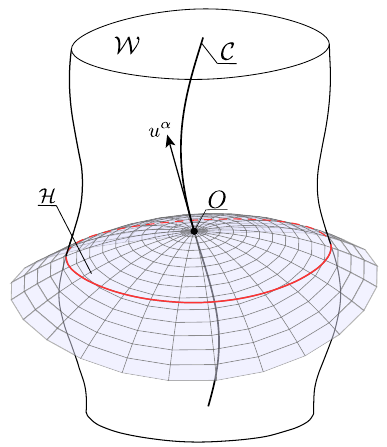}
\end{center}
\caption{ A schematic representation of the world tube.}
\label{fig2}
\end{figure}

	 Consider a gravitational field with a given metric in which the test body moves.
	 It is described by an energy-momentum tensor $T^{\alpha \beta}$ that satisfies
the following conditions:
	 \begin{itemize}
	 \item [(i)] the symmetry $T^{\alpha \beta} = T^{\beta \alpha}$;
	 \item [(ii)] the conservation law $\nabla_\beta T^{\alpha \beta}=0$.
	 \end{itemize}
The trajectory of the test body (see Fig. \ref{fig2}) on the manifold $\mathcal{M}^4$
is a world tube $\mathcal{W}={\rm supp}(T^{\alpha \beta})$
all spatial sections of which are compact.
Inside $\mathcal{W}$ we choose a supporting time-like curve
$
\mathcal{C}=\{ x^\alpha = y^\alpha (\tau)\},
$
where $\tau$ is an arbitrary (nonnatural) parameter.
 The tangent vector to this curve $u^{\alpha}=\dfrac{dy^\alpha}{d\tau}$ satisfies
 the condition $u^\alpha u_\alpha<0$.

Let $\mathcal{H}(\tau_*)$ be a space-like surface formed by combining all geodesics
$z^{\alpha}(s)$,
which, first, pass through the fixed point $O \in \mathcal{C}$ with
the parameter value $\tau=\tau_*$
and, second, have a tangent line orthogonal to $u^\alpha(\tau_*)$. Next, moving along the
supporting curve and varying the parameter $\tau$, we obtain a family of surfaces
$\mathcal{H}(\tau)$.

Assume that the spatial sections of the body are sufficiently small so that the following
condition is satisfied:
\begin{itemize}
 \item [(iii)] except at the points of the supporting curve $\mathcal{C}$, the geodesics
 $z^{\alpha}(s)$ constructed both for the fixed $\tau=\tau_*$ and for different values
 of $\tau$ do not intersect each other inside the tube $\mathcal{W}$.
\end{itemize}

Conditions (i)-(iii) imply that the evolution of the momentum of the test body, $p^\alpha$, and the tensor of the angular momentum, $S^{\alpha \beta}$, are described by the following system of equations (for details, see \cite{Harte2012}):
\begin{equation}
\label{eq_MP}
\begin{gathered}
\frac{Dp^\alpha}{d\tau}=-\frac{1}{2}R^{\alpha}_{ \ \lambda \mu \nu } u^{\lambda}S^{\mu \nu}, \\
\frac{DS^{\alpha \beta}}{d\tau} =p^{\alpha} u^{\beta} - u^{\alpha}p^{\beta},
\end{gathered}
\end{equation}	
where $R^{\alpha}_{ \ \mu \nu \lambda}$ is the curvature tensor and the following
covariant derivative has been introduced:
$$
\begin{gathered}
\frac{Dp^\alpha}{d\tau}=\frac{dp^\alpha}{d\tau} + \Gamma^\alpha_{\ \mu \nu} p^\mu u^\nu , \\
\frac{DS^{\alpha \beta}}{d\tau} = \frac{dS^{\alpha \beta}}{d\tau} + \Gamma^\alpha_{\ \mu \nu} S^{\mu \beta} u^\nu +
 \Gamma^\beta_{\ \mu \nu} S^{\alpha \mu} u^\nu,
\end{gathered}
$$	
where $\Gamma^\alpha_{\ \mu \nu}$ is the Christoffel symbol.

System \eqref{eq_MP} is usually referred to as the Mathisson\,---\,Papapetrou\,---\,Dixon (MPD) equations.
But in fact, Mathisson and Papapetrou obtained the equations in a somewhat different form. We consider their
results
in greater detail.

The problem of deriving the equations of motion for a rotating test body in general
relativity theory was first addressed by Mathisson (1937) \cite{Mathisson1937,Mathisson2010}.
He explicitly used the additional assumption
\begin{equation}
\label{eqMattison}
S^{\alpha \beta} u_\beta=0.
\end{equation}
 Mathisson obtained equations of motion equivalent to \eqref{eq_MP}, but taking
 condition \eqref{eqMattison} into account.
The work of Mathisson and
its further development are discussed in detail in Ref. \cite{Dixon2015}.

The equations of motion for a rotating body that are equivalent to \eqref{eq_MP} were
obtained by Papapetrou (1951)\cite{Papapetrou1951}.
He used a noncovariant definition of multipole moments. Later, his approach was refined in
Ref. \cite{TulczyjewTulczyjew}.
It should be noted that, in contrast to Mathisson, Papapetrou imposed no
additional condition, which allowed the search for such a condition
to be viewed as a separate topic.

The approaches of Mathisson and Papapetrou were simplified by Tulczyjew
\cite{Tulczyjew1959}(1959). He noted that the pole-dipole approximation reduces to an
energy-momentum tensor in the form
\begin{equation}
\label{T_Tulczyjew}
\begin{gathered}
T^{\alpha \beta }=\int\limits_\mathcal{C} d\tau \left[ \mu^{\alpha \beta} \delta^{(4)}(x - y(\tau)) - \mu^{\gamma \alpha \beta}\nabla_\gamma \delta^{(4)}(x - y(\tau)) \right], \\
\end{gathered}
\end{equation}
where $\delta^{(4)}(x - y(\tau)) $ is the four-dimensional Dirac delta-function.
Relations \eqref{T_Tulczyjew} can be regarded \cite{Steinhoff2011, Trautman2002} as
a relativistic generalization of the density in which the multipole
moments $\mu^{\alpha \beta}$ and $\mu^{\gamma \alpha \beta}$ are concentrated along the curve
$\mathcal{C}$.

Next, using some properties of integrals in the form \eqref{T_Tulczyjew}, Tulczyjew showed
that conditions (i) and (ii) lead to equations of motion that are equivalent to \eqref{eq_MP}.
The Tulczyjew method and the proof of the properties used by him
are discussed in detail in Ref. \cite{Ohashi2003}. Further development of the Tulczyjew method
taking into account multipole moments of higher order
is described in Ref. \cite{Steinhoff2010}.

A generalization of these results for an extended body, i.e., for a body with
an energy-momentum tensor that contains no Dirac delta-function, was made by Dixon
\cite{Dixon1964}(1964). He obtained equations of motion exactly in the form \eqref{eq_MP} and
gave a definition of the linear momentum and the angular momentum of the test body in curved space.
Later, Dixon obtained equations taking into account multipole moments of higher order
(for details, see \cite{Dixon1970, Dixon2015}). An explicit introduction of generalized
Killing fields and the generalized linear momentum by Harte
\cite{Harte2008} (2008) simplified the derivation of the equations of motion \eqref{eq_MP}.
In addition, Harte derived equations of motion in the Newtonian gravitation using
Killing fields in Ref. \cite{Harte2012}.

We note that under the pole-dipole approximation the equations do not explicitly depend on the
velocity distribution
in the test
body. If one takes into account multipole moments of higher order, this holds no longer
\cite{Dixon1974}.

\subsection{The Tulczyjew condition}
\label{sec_tul}
By assumption, $p^\alpha$ is time-like and $S^{\alpha \beta}$ is space-like. Hence,
\begin{equation}
\label{eq_mS}
m^2=-p_\alpha p^\alpha, \quad c^2=\frac{S_{\alpha \beta}S^{\alpha\beta}}{2} .
\end{equation}
Next, we let $m$ denote the mass of the body, and $c$,  
its spin magnitude. Also, we denote
$$
k = -p^\alpha u_\alpha.
$$
We point out a few important properties of the MPD equations
\eqref{eq_MP}.
\begin{itemize}
\item[---] If $S^{\alpha \beta}=0$, then this system reduces to an analysis of geodesics.
For them the linear momentum and the velocity are parallel to each other:
\begin{equation}
\label{eq_pu}
p^\alpha = m u^\alpha
\end{equation}
and $k=m$.
\item[---] If $S^{\alpha \beta}\neq 0$, then in the general case relation \eqref{eq_pu}
is not satisfied and $m\neq k$. Therefore, the linear momentum is not parallel to the velocity,
$p^\alpha \nparallel u^\alpha $, moreover, for them it is possible that $k=0$, i.e., they can
be orthogonal to each other (see Remark \ref{rem2}). However, if
condition \eqref{eq_pu} is fulfilled and $S^{\alpha \beta}\neq 0$, then it follows from
the MPD equations \cite{Kopeikin2011} that the angular momentum
    is transported in parallel along the trajectory. For example, this holds in de Sitter
    space \cite{ObukhovPuetzfeld2011}, which has a constant curvature.
\end{itemize}
As in the Newtonian gravitation \cite{Dixon1974}, system \eqref{eq_MP} turns out to be incomplete since
velocity $u^\alpha$ is not given as a function of $p^\alpha$ and $S^{\alpha \beta}$.
To determine velocity $u^\alpha$, we define the vector
$$
f^\alpha=S^{\alpha \beta}p_\beta.
$$
As noted by Tulczyjew \cite{Tulczyjew1959}, the additional condition \eqref{eqMattison}, which is
used by Mathisson, does not define a unique point even in flat space (for details, see
\cite{Costa2012}).
He suggested another condition, which we will call the {\it Tulczyjew condition} in what
follows:
\begin{equation}
\label{eq_TD}
f^\alpha=0.
\end{equation}
In Refs. \cite{Schattner1979a, Schattner1979b} it was proved that, if the linear and angular
momenta of a test body are to be determined as per Dixon's definition \cite{Dixon1976}, then condition \eqref{eq_TD} for
a extended body is satisfied by a unique point on the surface
$\mathcal{H}(\tau)$ under some additional, but completely natural, conditions imposed on
the energy-momentum tensor.

Condition \eqref{eq_TD} defines only three independent relations as
 $$
 f^\alpha p_\alpha = S^{\alpha \beta} p_\alpha p_\beta \equiv 0.
 $$
Note that Eqs. \eqref{eq_TD}  are a homogeneous linear system
in the components of the linear momentum $p_{\beta}$. Therefore, they have a nontrivial
solution if
\begin{equation}
\label{eq_T01}
\begin{gathered}
\det || S^{\alpha \beta} ||= \left( \frac{1}{8} \sqrt{-g} \varepsilon_{\alpha \beta \mu \nu} S^{\alpha \beta} S^{\mu \nu} \right)^2= (S^{0 1} S^{2 3} - S^{0 2} S^{1 3} + S^{0 3} S^{1 2})^2 =0,
\end{gathered}
\end{equation}
 whence it follows that $S^{\alpha \beta} $ satisfies the condition  \cite{Santos2020}
\begin{equation}
\label{eq_T02}
S^{\alpha \beta} S^{\mu \nu} - S^{\mu \beta} S^{\alpha \nu} + S^{\mu \alpha}S^{\beta \nu}=0.
\end{equation}
Another consequence of relation \eqref{eq_T01} is that the skew-symmetric matrix
$||S^{\alpha \beta}||$ has rank two.

The Tulczyjew condition implies that a reference system of the zero 3-momentum
$p^i=0$ has been chosen.
As a rule, such a supporting curve $\mathcal{C}$ is called the {\it  worldline}.

\begin{proposition}
\label{utv2}
If the Tulczyjew condition \eqref{eq_TD} is satisfied and $k\neq 0$, then system \eqref{eq_MP}
admits the additional integral
$m = {\rm const}$.
 \end{proposition}
\begin{proof}
Using system \eqref{eq_MP} and the asymmetry of the curvature tensor in the
first two indices, one can show that
$$
\frac{dm}{d\tau}=-\frac{1}{2m} R_{ \beta \alpha \mu \nu} p^{\beta} u^\alpha S^{\mu \nu}= \frac{p_\beta }{m k}\frac{DS^{\alpha \beta}}{d\tau}\frac{Dp_\alpha}{d\tau}.
$$
Next, we differentiate the vector $f^\alpha$ by virtue of system \eqref{eq_MP} and
set the resulting relation equal to zero
\begin{equation}
\label{eq_df}
\frac{df^\alpha}{d\tau}=p_\beta\frac{D S^{\alpha \beta}}{d\tau} + S^{\alpha \beta}\frac{D p_\beta}{d\tau} - \Gamma^\alpha_{\ \mu \nu} f^\mu u^\nu=0,
\end{equation}
where the last term is zero by virtue of \eqref{eq_TD}, whence
 we find that
$$
\frac{dm}{d\tau}=\frac{1}{m k}S^{\alpha \beta} \frac{Dp_\alpha}{d\tau} \frac{Dp_\beta}{d\tau}=0.
$$
\end{proof}
Relation \eqref{eq_df} is a linear system in the velocity components $u^\alpha$, which, using
\eqref{eq_TD}  and \eqref{eq_MP},
can be represented as
\begin{equation}
\label{eq_sys00}
\begin{gathered}
u^{\alpha} - \frac{1}{2 m^2} \tensor{\mathcal{D}}{^\alpha_\beta} u^\beta=\frac{k}{m^2} p^\alpha, \quad
\tensor{\mathcal{D}}{^\alpha_\beta} = S^{\alpha \gamma} R_{\gamma \beta \mu \nu} S^{\mu \nu}.
\end{gathered}
\end{equation}
To solve this system, we use the following statement \cite{Santos2020}.
\begin{proposition}
\label{utv20}
For the matrix ${\bf D}= || \tensor{\mathcal{D}}{^\alpha_\beta} ||$,
the following holds true:
\begin{equation}
\label{eq_D2}
{\bf D}^2 = d {\bf D}, \quad d = \dfrac{1}{2}{\rm tr}({\bf D})=-\frac{1}{2} R_{\alpha \beta \mu \nu} S^{\alpha \beta }S^{\mu \nu}.
\end{equation}
 \end{proposition}
\begin{proof}
We rewrite the elements of the matrix ${\bf D}$ as follows:
$$
\tensor{\mathcal{D}}{^\alpha_\beta} = S^{\alpha \gamma}Q_{\gamma \beta }, \quad Q_{\gamma \beta} = R_{\gamma \beta \mu \nu} S^{\mu \nu}.
$$
For the elements of the matrix ${\bf D}^2$, using the property \eqref{eq_T02} we obtain
$$
\begin{gathered}
\tensor{\mathcal{D}}{^\alpha_\mu} \tensor{\mathcal{D}}{^\mu_\gamma} = S^{\alpha \beta} S^{\mu \nu} Q_{\beta \mu} Q_{\nu \gamma} = (S^{\mu \beta} S^{\alpha \nu} - S^{\mu \alpha} S^{\beta \nu})Q_{\beta \mu} Q_{\nu \gamma} = \tensor{\mathcal{D}}{^\mu_\mu} \tensor{\mathcal{D}}{^\alpha_\gamma} -
\tensor{\mathcal{D}}{^\alpha_\beta} \tensor{\mathcal{D}}{^\beta_\gamma},
\end{gathered}
$$
where account is taken of the fact that $Q_{\beta \mu}= - Q_{\mu \beta}$.
As we see, it follows from this relation that
$$
{\bf D}^2 = {\rm tr}({\bf D}) {\bf D} - {\bf D}^2,
$$
which leads us to \eqref{eq_D2}.
\end{proof}
Using the Neumann series and the property \eqref{eq_D2}, we obtain
$$
\left( 1 - \frac{{\bf D}}{2m^2} \right)^{-1} = \sum_{n=0}^{\infty} \frac{{\bf D}^n}{(2m^2)^n} = 1 - \frac{{\bf D}}{d - 2m^2}.
$$
Multiplying system \eqref{eq_sys00} by this matrix, we find that for $d \neq 2 m^2 $
the velocity is given by \begin{equation}
\label{eq_vel}
\begin{gathered}
u^\alpha= \frac{k}{m^2}\left( p^{\alpha} + \dfrac{\tensor{\mathcal{D}}{^\alpha_\beta} p^{\beta}}{ 2 m^2 - d} \right).
\end{gathered}
\end{equation}
In order to completely define the velocity, we need to define the parameter $\tau$.
We will stick to the parametrization proposed in Ref. \cite{Ehlers1977}(see also
\cite{LukesGerakopoulos2017}), according to which this
parameter is chosen in such a way that
\begin{equation}
\label{eq_params}
\frac{k}{m}=1.
\end{equation}
Velocity $u^\alpha$ is a time-like vector. Hence, it follows from \eqref{eq_vel},
\eqref{eq_params} and \eqref{eq_TD} that
\begin{equation}
\label{eq_ner1}
u_\alpha u^\alpha = \frac{\tensor{\mathcal{D}}{^\alpha_\beta} \mathcal{D}_{\alpha \gamma}p^\beta p^\gamma}{(d - 2 m^2)^2} - m^2 <0.
\end{equation}
As is well known, this condition can generally be violated (for examples, see
\cite{Tod1976, Hojman1975}). If this condition is not satisfied, then the pole-dipole
approximation is violated~\cite{Semerak1, HojmanHojman} and one needs to consider
multipole moments of higher order. Also, Ref. \cite{Deriglazov2020} shows that, in an arbitrary gravitational field, the MPD equations can have solutions for which the acceleration increases with the velocity and diverges in the ultrarelativistic limit. A modification of these equations that exhibits correct behavior in the ultrarelativistic limit is discussed in Ref. \cite{GuzmnRamrez2022}.

The MPD equations \eqref{eq_MP} for velocity \eqref{eq_vel} and
\eqref{eq_params}
reduce to the following system:
\begin{equation}
\label{eq_sys}
\begin{aligned}
\frac{DS^{\alpha \beta}}{d\tau} &= \frac{p^\alpha p^\beta - p^\beta p^\alpha}{m} + \frac{p^\alpha\tensor{\mathcal{D} }{^\beta_\mu} - p^\beta\tensor{\mathcal{D} }{^\alpha_\mu}}{m(d - 2m^2)}p^\mu ,\\
\frac{Dp_\alpha}{d\tau} &= -\frac{1}{2m}R_{\alpha\beta \mu \nu}\left[ p^\beta + \frac{\tensor{\mathcal{D}}{^\beta_\lambda} p^\lambda}{ d - 2m^2} \right]S^{\mu \nu},\\
\frac{d y^\alpha}{d\tau} &= \frac{p^\alpha}{m} + \frac{\tensor{\mathcal{D}}{^\alpha_\beta} p^\beta}{m(d - 2 m^2 )}.
\end{aligned}
\end{equation}
If the metric $g_{\alpha \beta}$ has a Killing field $\xi^\alpha$
then the equations of motion \eqref{eq_sys} have the additional integral \cite{Dixon1976}:
\begin{equation}
\label{eq_F}
F(\xi ) = p_\alpha \xi^\alpha + \frac{1}{2}S^{\alpha \beta}\nabla_\alpha \xi_\beta.
\end{equation}

Of special note is that the equations for velocity \eqref{eq_sys00} explicitly take
the Tulczyjew condition into account. If one finds velocity $u^\alpha$ from the solution
to the initial equation \eqref{eq_df}, then the MPD equations
will have the additional integrals $f^\alpha={\rm const}$. This general solution is
rather cumbersome, but for $f^\alpha=0$ it reduces to the solution \eqref{eq_vel}.
Thus, the analysis of system \eqref{eq_sys} should be carried out exclusively on its
invariant manifold given by \eqref{eq_TD}. Specifically, the stability of the partial
solutions should be investigated only with respect to perturbations lying on this manifold.

\begin{rem}
\label{rem2}
System \eqref{eq_sys} contains a singularity if $d - 2 m^2=0$. In this case,
one should seek a different solution of system \eqref{eq_sys00}.
For this solution, it may turn out that $k=0$ and mass $m$ will no longer remain unchanged. We give the following
example.

We first note that, for the curved space-time and the nonzero
angular momentum $S^{\alpha \beta } \neq 0$,
the condition $k=0$ does not lead to null geodesics since, according to \eqref{eq_sys00}, we
have
$$
u_\alpha u^\alpha = \frac{1}{2m^2} \tensor{\mathcal{D}}{_\alpha_\beta} u^\alpha u^\beta \neq 0.
$$

Suppose that the initial conditions lie on the surface
\begin{equation}
\label{eqK}
N= p^\alpha p_\alpha - \frac{1}{4} R_{\alpha \beta \mu \nu} S^{\alpha \beta }S^{\mu \nu} = 0.
\end{equation}
Then, taking $k=0$ into account, Eq. \eqref{eq_sys00} has the form
$$
\tensor{\mathcal{A}}{^\alpha_\beta} u^\beta =0, \quad \tensor{\mathcal{A}}{^\alpha_\beta} = \tensor{\delta}{^\alpha_\beta} - \frac{\tensor{\mathcal{D}}{^\alpha_\beta}}{d}.
$$
If one considers the Schwarzschild metric, then ${\rm rank} || \tensor{\mathcal{A}}{^\alpha_\beta} ||=2$
and hence the general solution can be represented as
\begin{equation}
\label{eq_u0}
u^\alpha=j u_o^\alpha + l u_*^\alpha,
\end{equation}
where $u_o^\alpha$ and $u_*^\alpha$ are two linearly independent partial solutions and
$j,l={\rm const}$.
To define one of these constants, we require that the
MPD equations \eqref{eq_MP} admit an additional integral $K$ in the case of a
Schwarzschild metric. To
do so we differentiate it by virtue of this system and set it equal to zero
$$
\frac{dN}{d\tau}=2R_{\alpha \beta \mu \nu} p^\alpha u^\beta S^{\mu \nu} + \frac{1}{4}S^{\alpha \beta} S^{\mu \nu} u^\lambda \nabla_\lambda R_{\alpha \beta \mu \nu}=0.
$$
Substituting velocity \eqref{eq_u0} into this equation, we obtain a linear equation
from which we express one of the constants. The other constant can be found from the
normalization condition, for example, $u_\alpha u^\alpha=-1$. Then, substituting the
velocity, we obtain MPD equations that are different from
\eqref{eq_sys}. By direct calculations one can show that, in the case of
a Schwarzschild metric, mass $m$ does not remain unchanged for this system.
\end{rem}

\subsection{Hamiltonian form}
The resulting system \eqref{eq_sys} governs the evolution of $S^{\alpha \beta}$, $p_\alpha$
and $y^\alpha$.
We define the Poisson brackets between these variables as follows:
\begin{equation}
\begin{gathered}
\label{eq_J}
\{ S^{\alpha \beta}, S^{\mu \nu} \} = g^{\alpha\mu}S^{\beta \nu} - g^{\alpha\nu}S^{\beta \mu} + g^{\beta\nu}S^{\alpha \mu} - g^{\beta\mu}S^{\alpha \nu}, \\
\{ S^{\alpha \beta}, p_\mu \} = -\Gamma^{\alpha}_{\lambda \mu}S^{\lambda \beta} - \Gamma^{\beta}_{\lambda \mu}S^{ \alpha \lambda}, \\
\{ y^{\alpha}, p_\beta \}=\delta^\alpha_\beta, \quad \{ p_\alpha, p_\beta \} = -\frac{1}{2}R_{\alpha \beta \mu \nu}S^{\mu \nu},
\end{gathered}
\end{equation}
where only nonzero Poisson brackets are presented. This Poisson bracket arose first in
the work of Souriau \cite{Souriau1970} (1970) (for a detailed discussion of his work, see
also \cite{Damour2024}).
Later, the bracket \eqref{eq_J} was independently found in Refs. \cite{Feldman1980, Kunzle1972, Khriplovich1989, Kumar2015}.

The equations of motion \eqref{eq_sys} on the invariant manifold \eqref{eq_TD}
 can be represented as
\begin{equation}
\label{eq_JJ}
\begin{aligned}
\frac{d S^{\alpha \beta}}{d\tau}&=\{ S^{\alpha \beta}, S^{\mu \nu} \} \frac{\partial H}{\partial S^{\mu \nu}} + \{ S^{\alpha \beta}, p_\mu \} \frac{\partial H}{\partial p_\mu}, \\
\frac{d p_\alpha}{d\tau}&= \{ p_\alpha, S^{\mu \nu} \} \frac{\partial H}{\partial S^{\mu \nu}} + \{ p_\alpha, p_\mu \} \frac{\partial H}{\partial p_\mu} - \frac{\partial H}{\partial y^{\alpha}} ,\\
\frac{d y^\alpha}{d\tau} &= \frac{\partial H}{\partial p_\alpha},
\end{aligned}
\end{equation}
where the Hamiltonian $H$ is given by\footnote{A passage from the Hamiltonian system \eqref{eq_JJ}, \eqref{eq_Hamilton} with the Dirac
constraints \eqref{eq_TD} to an equivalent (degenerate) Lagrangian system was examined in detail
in Ref. \cite{Ramirez2015}.}
\begin{equation}
\label{eq_Hamilton}
H=\frac{ p^\alpha p_\alpha}{2 m} - \dfrac{\mathcal{D}^{\alpha \beta} p_{\alpha}p_\beta}{m(d + 2 p_\alpha p^\alpha)} .
\end{equation}
On the invariant relation \eqref{eq_TD} the Hamiltonian takes the fixed value $H=-\dfrac{m}{2}$.

The Poisson tensor defined by the bracket \eqref{eq_J} has rank 12, thus giving rise to
two Casimir functions:
$$
\begin{gathered}
C_\star=\frac{1}{8}\sqrt{-g} \varepsilon_{\alpha \beta \mu \nu} S^{\alpha\beta}S^{\mu \nu}, \quad
C_\circ=\frac{1}{2}S^{\alpha\beta}S_{\alpha \beta}.
\end{gathered}
$$
According to \eqref{eq_T01}, one Casimir function has zero value, and the other takes,
according to \eqref{eq_mS}, positive values, so that the trajectories we are interested in lie
on the family of symplectic leaves
$$
\mathcal{M}^{12}_c=\{ C_\star(S^{\alpha \beta}) = 0, C_\circ( S^{\alpha \beta})= c^2 >0 \}.
$$

The Poisson bracket \eqref{eq_J} can be simplified if we transform from the coordinate
basis in the tangent plane $T\mathcal{M}^4$
to the orthonormal basis composed of the tetrad $\{ \tensor{e}{_A^{\alpha}} \}$, four
vector fields satisfying the relations
$
g_{\alpha \beta} \tensor{e}{_A^{\alpha}}\tensor{e}{_B^{\beta}}= \eta_{AB},
$
where $\eta_{AB}={\rm diag}(-1, 1,1,1)$. Here and in the sequel the indices
labeled by capital Latin letters will take values from 0 to 3 and denote the
sequence number of the tetrad.

Let $P_\alpha$ denote the new linear momentum, which we define as follows
(for details, see \cite{Tauber1988, Feldman1980}):
\begin{equation}
\label{eq_Sp}
p_\alpha = P_\alpha + \frac{1}{2}\omega_{\alpha AB} S^{AB}, \quad
S^{\alpha \beta} = S^{AB} \tensor{e}{_A^\alpha} \tensor{e}{_B^\beta},
\end{equation}
where the Ricci coefficients of rotation have been introduced
\begin{equation}
\label{eq_omega}
\omega_{\alpha AB}=g_{\mu \nu}\tensor{e}{_A^{\mu}} \nabla_\alpha \tensor{e}{_B^\nu}.
\end{equation}
Hence, the Poisson bracket between the new linear momenta will be equal to zero:
$\{ P_\alpha, P_\beta \}=0$.

Define two three-dimensional vectors
$$
\begin{gathered}
\vec{L}=(L_1, L_2, L_3) = (S^{01}, S^{02}, S^{03}), \\
\vec{M}=(M_1, M_2, M_3) = (S^{23}, S^{31}, S^{12}).
\end{gathered}
$$
Thus, these vectors parameterize the components of the matrix $S^{AB}$ in the
orthonormal basis
\begin{equation}
\label{eq_SAB}
S^{AB}= \left(
\begin{array}{cccc}
0 & \vec{M}\\
-\vec{M}^{T} & \hat{L}
\end{array}
\right).
\end{equation}
 Such a parameterization of the matrix $S^{AB}$ is sufficiently convenient and
 has already been used in Refs. \cite{Ramond2025, Ramond2022}.
 In the new variables the nonzero Poisson brackets \eqref{eq_J} become
\begin{equation}
\label{eq_J2}
\begin{gathered}
\{ L_i, L_j\} = \hat{L}_{ij}, \quad
\{ L_i, M_j\} = \hat{M}_{ij}, \quad
\{ M_i, M_j\} = - \hat{L}_{ij} \\
\{ y^{\alpha}, P_\beta \} =\tensor{\delta}{^\alpha_\beta},
\end{gathered}
\end{equation}
As we see, the commutation relations between $\vec{L}$ and $\vec{M}$ form the algebra
$so(3,1)$, for which the preceding Casimir functions take the form
\begin{equation}
\label{eq_C2}
C_\star=(\vec{L}, \vec{M}), \quad
C_\circ=(\vec{L}, \vec{L}) - (\vec{M}, \vec{M}).
\end{equation}
Finally, the equations of motion can be represented in the following vector form
(for details, see Appendix \ref{appendicesA}):
\begin{equation}
\label{eq_HamSys}
\begin{aligned}
\frac{d \vec{L}}{d \tau} &=\frac{\partial H}{\partial\vec{L}}\times \vec{L} + \frac{\partial H}{\partial\vec{M}} \times \vec{M}, \\
\frac{d \vec{M}}{d \tau} &= \vec{M}\times\frac{\partial H}{\partial\vec{L}} + \frac{\partial H}{\partial\vec{M}} \times \vec{L}, \\
\frac{d P_\alpha}{d \tau}&= - \frac{\partial H}{\partial y^{\alpha}} ,\quad
\frac{d y^\alpha}{d \tau} = \frac{\partial H}{\partial P_\alpha},
\end{aligned}
\end{equation}
where $H=H(\vec{L}, \vec{M}, P_\alpha, y^\alpha)$ is the Hamiltonian in the new variables.
Direct substitution of \eqref{eq_Sp} and \eqref{eq_SAB} into \eqref{eq_Hamilton} and \eqref{eq_TD} gives explicit relations for $H$ and
$f^\alpha$ in the new variables. For an arbitrary metric, these expressions are
rather cumbersome and cannot be simplified to a compact form, and so they are not presented here.
However, in specific cases, for example, in a Schwarzschild
metric or in a Kerr metric, the expressions $H$ and $f^\alpha$ can be obtained
using any system of analytic calculations (for example, Maple or Mathematica).

Thus, in the general case system \eqref{eq_HamSys} is a Hamiltonian system with six
degrees of freedom \cite{Witzany2019}. The main advantage of the chosen variables is that
the Poisson bracket in them is the simplest and the
equations for $\vec{M}$ and $\vec{L}$ are similar to systems arising in rigid body dynamics \cite{BorisovMamaev2018,BorisovMamaevNonEuclideanSpaces, Tsiganov2024}.
All this allows one to apply efficient reduction methods to system \eqref{eq_HamSys} in
the case of additional integrals of the form \eqref{eq_F}. Below this will be illustrated
by a Schwarzschild metric. Previously, such reduction methods have proved to be efficient
in the problem of two bodies in spaces of
constant curvature \cite{Arathoon2023, BorisovMamaevBizyaev}.

\section{The Schwarzschild space-time}
If we denote the coordinates by $(x^\alpha)=(t, r, \theta, \varphi)$, then the interval for
a Schwarzschild metric can be represented as
\begin{equation}
\label{eq_sw}
ds^2=- a(r)dt^2 + \frac{dr^2}{a(r)} + r^2(d\theta^2 + \sin^2 \theta d\varphi^2), \
a(r)=1 - \frac{2\mu}{r}.
\end{equation}
This interval describes a spherically symmetric black hole with the event horizon
$
\mathcal{S}_h=\{ (x^\alpha) \ | \ r=2\mu\},
$
where the parameter $\mu$ has the physical meaning of the mass of a black hole.
Below we will consider the motion of the test body only on the outside of the event horizon.

The nonzero components of the curvature tensor in the basis attached to the coordinates
have the form
$$
\begin{gathered}
R_{trtr}=-\frac{2\mu}{r^3}, \quad R_{t \theta t \theta}= a(r)\frac{\mu}{r}, \quad R_{t \varphi t \varphi} = a(r)\frac{\mu}{r}\sin^2\theta, \\
R_{r \theta r \theta}=- \frac{\mu}{ra(r)}, \quad
R_{r \varphi r \varphi} = - \frac{\mu}{ra(r)}\sin^2\theta, \quad
R_{\theta \varphi \theta \varphi} = 2\mu r \sin^2 \theta,
\end{gathered}
$$
and the other nonzero components can be obtained from the antisymmetry properties
$R_{\alpha \beta \mu \nu}= - R_{\alpha \beta \nu \mu}=-R_{\beta \alpha \mu \nu}= R_{\beta \alpha \nu \mu }$.

Since the metric is diagonal, orthonormalized vectors can be obtained by normalizing
the chosen basis
$$
e_0=\frac{1}{\sqrt{a(r)}}\frac{\partial }{\partial t}, \ e_1 = \sqrt{a(r)} \frac{\partial }{\partial r}, \
e_2=\frac{1}{r}\frac{\partial}{ \partial \theta}, \
e_3=\frac{1}{r \sin \theta} \frac{\partial }{\partial \varphi}.
$$

Next, we will need relations for the coefficients $\omega_{\alpha AB}$. Straightforward calculation
from \eqref{eq_omega} shows that only eight coefficients do not vanish:
$$
\begin{gathered}
\omega_{t01}=-\frac{\mu}{r^2}, \ \omega_{\theta 12} = -\sqrt{a(r)}, \
 \omega_{\varphi 1 3}=-\sqrt{a(r)}\sin \theta, \
\omega_{\varphi 2 3}=-\cos\theta,
\end{gathered}
$$
and four other coefficients can be obtained from the antisymmetry property
$\omega_{\alpha AB}=-\omega_{\alpha B A}$.

In the case of a Schwarzschild metric, system \eqref{eq_sys} has the discrete symmetry:
\begin{equation}
\label{eq_Sym}
\begin{gathered}
\mathcal{I}: \ \theta \to \pi-\theta, \ p_\theta \to -p_\theta, \ S^{r \theta} \to -S^{r \theta}, \ S^{t \theta} \to -S^{t \theta}, \ S^{\theta \varphi} \to - S^{\theta \varphi},
\end{gathered}
\end{equation}
where the other variables remain unchanged. As is well known, the fixed symmetry points
form an invariant submanifold which for
$\mathcal{I}$ consists of trajectories lying in the equatorial plane
\begin{equation}
\label{eq_th_p_S}
 {\rm Fix } \ \mathcal{I}: \ \theta=\frac{\pi}{2}, \quad p_\theta = 0, \quad S^{r \theta}=S^{t \theta}=S^{\theta \varphi}=0.
\end{equation}

Next, the explicit coordinates of the tensors in the coordinate basis will be denoted by
$t, r, \theta, \varphi$, and the components of the three-dimensional vectors $\vec{L}$ and
$\vec{M}$ will be labeled by $1,2,3$ to emphasize their difference.

\subsection{ Additional integrals}
The metric \eqref{eq_sw} has four Killing fields:
\begin{equation}
\label{eqKilling}
\begin{gathered}
 \frac{\partial}{\partial t}, \ -\sin \varphi \frac{\partial}{\partial \theta} - \frac{\cos \varphi}{\tan \theta} \frac{\partial}{\partial \varphi}, \
 \cos \varphi \frac{\partial}{\partial \theta} - \frac{\sin \varphi}{\tan \theta} \frac{\partial}{\partial \varphi}, \
 \frac{\partial}{\partial \varphi}.
\end{gathered}
\end{equation}
As a result, system \eqref{eq_HamSys} admits four additional integrals, which, according
to \eqref{eq_F}, can be represented as
\begin{equation}
\label{eq_Int}
\begin{gathered}
\mathcal{E} = -P_t, \quad \vec{Q}=(Q_1, Q_2, P_\varphi), \\
Q_1= \frac{\cos \varphi}{\sin \theta}(L_1 - P_\varphi \cos \theta ) - P_\theta\sin\varphi, \\
Q_2= \frac{\sin \varphi}{\sin \theta}(L_1 - P_\varphi \cos \theta ) + P_\theta \cos\varphi.
\end{gathered}
\end{equation}
where $\mathcal{E}$ is the energy of the test body.
The last three Killing fields in \eqref{eqKilling} are the generators of the rotation group $SO(3)$,
and so the components of $\vec{Q}$ can be interpreted \cite{Suzuki1997}  as components of the total angular momentum (spin+orbital) of the test body.  

The Poisson bracket \eqref{eq_J2} of the integral $\mathcal{E}$ with the components of the
vector $\vec{Q}$ is zero.
Between the components of $\vec{Q}$ this bracket has the form
$$
\{ Q_i, Q_j\} = \hat{Q}_{ij}.
$$
As we see, if $\vec{Q} \neq 0$, then the set of integrals \eqref{eq_Int} is noninvolutive
and forms the algebra $so(3)$. From the components of $\vec{Q}$ one can combine two
involutive integrals $P_\varphi$ and
 \begin{equation}
 \begin{gathered}
\label{eq_001}
F= (\vec{Q}, \vec{Q}) = P_\theta^2 + P_\varphi^2 + \frac{(L_1 - P_\varphi \cos \theta)^2}{\sin^2 \theta} = P_\theta^2 + L_1^2 + \frac{(L_1 \cos \theta - P_\varphi)^2}{\sin^2 \theta} = \ell \geqslant 0.
\end{gathered}
\end{equation}
The integral $F$ is the total angular momentum of the test body. Throughout the remainder of the paper, we will denote the magnitude of the total angular momentum by $\ell$. 
Using the general result on order reduction for
Hamiltonian systems (see, e.g., \cite{ArnoldKozlov}), we arrive at the following conclusion.

\begin{proposition}
The Hamiltonian system \eqref{eq_HamSys} for a Schwarzschild metric admits the reduction
 \begin{itemize}
 \item [---] by three degrees of freedom if $\vec{Q} \neq 0$ since
 the integrals $Q_3$, $\mathcal{E}$ and $F$ are in involution with each other;
 \item [---] by four degrees of freedom if $\vec{Q} = 0$  since the integrals $Q_1$, $Q_2$, $Q_3$ and $\mathcal{E}$ are in involution with each other.
 \end{itemize}
\end{proposition}

Next, we consider the reduction in these two cases separately.

\subsection{Reduction in the case $\vec{Q}\neq 0$}
\label{sectionRedQneq0}
The metric \eqref{eq_sw} does not explicitly depend on $\varphi$ and $t$. Therefore,
reduction by two degrees of freedom has, in fact, been performed in system \eqref{eq_HamSys}
as the system governing the evolution of $\vec{L}$, $\vec{M}$, $P_r$, $r$, $P_\theta$,
and $\theta$ decouples.
To perform further reduction of the order of the system, we write the symmetry field
generated by the
integral $F$:
\begin{equation}
\label{eq_symu}
\begin{gathered}
\frac{P_{\varphi} \cos \theta - L_1}{\sin^2 \theta} \left[ L_3 \frac{\partial }{\partial L_2} - L_2 \frac{\partial }{\partial L_3} + M_3 \frac{\partial}{\partial M_2} - M_2 \frac{\partial }{\partial M_3} - \frac{L_1 \cos \theta - P_{\varphi}}{\sin\theta} \frac{\partial}{\partial P_\theta} \right] + P_\theta \frac{\partial }{\partial \theta}.
\end{gathered}
\end{equation}
As is well known (for details see, e.g., \cite{ArnoldKozlov}), to perform the reduction
one needs to choose as new variables the integrals of this symmetry field, i.e., the
variables in which the Poisson bracket \eqref{eq_J2} with integral $F$ is equal to zero.
However, arbitrariness in the choice of reduced variables arises since
an arbitrary function of these integrals is also an integral.

In this case, the vector field \eqref{eq_symu} is fairly simple and it is easy to find
its independent additional integrals:
$$
\begin{gathered}
 L_1, \ \sqrt{L_2^2 + L_3^2}, \ (\vec{M}, \vec{L}), \ M_2L_3 - M_3L_2, \\
L_3 P_\theta - \frac{L_2}{\sin \theta}(L_1 \cos \theta - P_\varphi), \ L_2 P_\theta + \frac{L_3}{\sin \theta}(L_1 \cos \theta - P_\varphi).
\end{gathered}
$$
We make use of the arbitrariness described above to simplify the Poisson bracket
for the value of the Casimir function, $C_\star=0$.
By straightforward calculations one can verify that it is convenient to choose
$P_r$, $r$ and \begin{equation}
\label{eq_EZ0}
\begin{gathered}
\vec{E}= \left( L_1, \frac{L_3 P_\theta}{\sqrt{L_2^2 + L_3^2}} + L_2\frac{L_1 \cos \theta - P_\varphi}{\sin \theta \sqrt{L_2^2 + L_3^2}},
 -\frac{L_2 P_\theta}{\sqrt{L_2^2 + L_3^2}} + L_3\frac{L_1 \cos \theta - P_\varphi}{\sin \theta \sqrt{L_2^2 + L_3^2}} \right), \\
\vec{Z}=\left(-\frac{M_1 \sqrt{(\vec{L}, \vec{L})} }{\sqrt{L_2^2 + L_3^2}}, \sqrt{(\vec{L}, \vec{L})}, \frac{M_2L_3 - M_3L_2}{\sqrt{L_2^2 + L_3^2}} \right),
\end{gathered}
\end{equation}
as new variables since in these variables the nonzero Poisson brackets on the level set of the
integral $C_\star=0$ have a fairly simple form (for details, see Appendix \ref{appendicesA}):
\begin{subequations}
\label{SP}
\begin{equation}
\label{SPa}
\begin{aligned}
\{ E_1, E_2 \} &= E_3, \quad \{ E_1, E_3 \} = -E_2, \quad \{ E_2, E_3 \} = E_1, \\
\{ Z_1, Z_2 \} &= Z_3, \quad \{ Z_1, Z_3 \} = Z_2, \quad \quad \{ Z_2, Z_3 \} = Z_1,
\end{aligned}
\end{equation}
\begin{equation}
    \{ r, P_r\} = 1.
\end{equation}
\end{subequations}
Note that the Poisson bracket \eqref{SPa} corresponds to expansion in the direct sum
$so(3)\otimes so(2,1)$.

The rank of the matrix defined by the Poisson bracket \eqref{SP} is six. Hence, it has
two Casimir functions. One of them is $C_\circ$, which is defined by \eqref{eq_C2},
and the other is the additional integral \eqref{eq_001}.
In the coordinates \eqref{eq_EZ0} they have the form
\begin{equation}
\label{eq_Co_F}
C_\circ = -Z_1^2 + Z_2^2 - Z_3^2, \quad
F = E_1^2 + E_2^2 + E_3^2. \quad
\end{equation}
In the new variables, the Hamiltonian \eqref{eq_Hamilton} becomes
\begin{equation}
\label{eq_Hi}
\begin{aligned}
H = & \frac{a(r)}{m} \left[ \frac{1}{2} + \frac{2\Theta^2 + Z_2^2 - E_1^2 }{V} \right] P_r^2 +
\frac{\Gamma + E_2}{mr V}\left[E_1 E_3 \left( \frac{Z_1^2}{Z_2^2} - 1\right) - Z_3 E_2\frac{Z_1}{Z_2} \right] P_r
 + \\ &\frac{Z_2^2 - E_1^2}{m r^3 V} \big[ r\Gamma^2 + \mu Z_1^2 \big] -
  \frac{( \mu \Theta - r^2 \mathcal{E})^2}{m r^3(r - 2\mu)}\left[ \frac{1}{2} - \frac{Z_1^2 + Z_3^2- 3\Theta^2}{V}\right] + \\ &
\frac{\sqrt{Z_2^2 - E_1^2}}{m r V}\big[ Z_3(\Theta - 2 r \mathcal{E}) P_r - Z_1Z_2 \mathcal{E} \big] +
\frac{E_3^2 + \Gamma^2}{m r^2}\left[ \frac{1}{2} - \frac{2E_1^2}{V} \right] - \\& \frac{ \mu \Theta - r^2 \mathcal{E} }{ m r^{2} V\sqrt{r(r - 2 \mu)}} \big[E_1 E_3 Z_3 + Z_1Z_2 E_2\big] - \frac{1}{m r^2V} \left[ \frac{Z_1}{Z_2}E_1E_3 + Z_3\Gamma \right]^2,
\end{aligned}
\end{equation}
where the following notation has been introduced:
$$
\begin{gathered}
V=\frac{m^2}{\mu}r^3 - 2Z_1^2 - Z_2^2 + Z_3^2 + 3 E_1^2\left(1 + \frac{Z_1^2 }{Z_2^2}\right), \\
\Gamma= \sqrt{a(r)(Z_2^2 - E_1^2)} - E_2, \quad
\Theta= \frac{Z_1}{Z_2}\sqrt{Z_2^2 - E_1^2}.
\end{gathered}
$$
Thus, system \eqref{eq_HamSys} reduces to an analysis of the following Hamiltonian system
with three degrees of freedom:
\begin{equation}
\label{eq_EZ}
\begin{gathered}
\frac{d \vec{E}}{d\tau} = \vec{E} \times \frac{\partial H}{\partial \vec{E}}, \quad
\frac{d \vec{Z}}{d\tau} = ({\bf J} \vec{Z}) \times \frac{\partial H}{\partial \vec{Z}}, \\
\frac{d P_r}{d\tau}= - \frac{\partial H}{\partial r} ,\quad
\frac{d r}{d\tau} = \frac{\partial H}{\partial P_r},
\end{gathered}
\end{equation}
where ${\bf J}={\rm diag}(1,-1,1)$. Note that, in deriving the equations of
motion \eqref{eq_EZ} and in the other relations in this section in which the partial
derivative of the Hamiltonian appears, we need to explicitly substitute the following
equation into the function $V$:
\begin{equation}
\label{eq_mmm}
m^2 = -a(r) P_r^2 + \frac{(\mu \Theta - \mathcal{E} r^2)^2}{r^3(r - 2 \mu)} - \frac{E_3^2 + \Gamma^2}{r^2}.
\end{equation}
Since $C_\circ>0$, the trajectories of system \eqref{eq_EZ} lie on the family of
symplectic leaves
\begin{equation}
\label{eq_ES6c}
\mathcal{S}^6_c= \{ (\vec{E}, \vec{Z}, r, P_r) \ | \ C_\circ(\vec{Z})=c^2, \ F(\vec{E})=\ell \}.
\end{equation}
Note that the Hamiltonian \eqref{eq_Hi} contains the term $\dfrac{Z_1}{Z_2}$. Therefore, it
has a singularity if $Z_2 = 0$, while on the leaves $\mathcal{S}^6_c$ the coordinates
of the vector $\vec{Z}$ lie on a two-sheet hyperboloid whose cavities are separated
by the plane $Z_2=0$, so that in this case for the trajectories we have $Z_2 \neq 0$.
In addition, according to \eqref{eq_EZ0}, if $\vec{L}\neq 0$ then $Z_2>0$.

We identify the singularities of Eqs. \eqref{eq_EZ} which can be reached by the trajectories
of this system. First, they include the event horizon $\mathcal{S}_h$, which corresponds to
the value of the radial coordinate $r=2 \mu$. Second, they include the case $d = 2 m^2 $
pointed out in Sec. \eqref{sec_tul}, which for the reduced
system is equivalent to $V=0$. The above-mentioned two singularities were present in the initial
system \eqref{eq_HamSys}. In addition to them, the equations of motion
have a singularity after reduction if $E_1=\pm Z_2$, which in the initial system reduces
to the case $L_2=0$, $L_3=0$. This singularity requires a separate analysis.

System \eqref{eq_EZ} has an invariant manifold defined by two relations from the
Tulczyjew condition:
\begin{equation}
\label{eq_Tred}
\begin{gathered}
f^r=-\Theta P_r - \sqrt{Z_2^2 - E_1^2}\frac{Z_3}{r} + \frac{E_2 Z_2 Z_3 - E_1 Z_1 E_3 }{Z_2\sqrt{r(r - 2 \mu)}} = 0, \\
f^\theta = \mathcal{E} \Theta + \frac{\mu \Theta^2}{r^2} - (\Gamma + E_2)\frac{\Gamma}{r}=0.
\end{gathered}
\end{equation}
On this invariant manifold, system \eqref{eq_EZ} describes a four-dimensional flow on
which the desired trajectories lie. This flow is given by a Hamiltonian system with two
degrees of freedom since one degree of freedom can be eliminated using Eqs. \eqref{eq_Tred}.
Hence, for it to be integrable by quadratures, we need an additional
integral. In the general case, this integral is absent \cite{Suzuki1997} and the problem
reduces to an analysis of a two-dimensional Poincar\'{e} map.

We now consider how one can retrieve the
values of the coordinates $\theta$, $\varphi$ and $t$ from the known solution of the reduced
system \eqref{eq_EZ}.
Since the Schwarzschild space-time is spherically symmetric, for each trajectory
there exists a coordinate system \cite{Suzuki1997} in which
\begin{equation}
\label{eq_Q12}
Q_1=0, \quad Q_2=0,
\end{equation}
 i.e., the total angular momentum is directed along the z-axis. 
This gives
\begin{equation}
\label{eqth}
P_\theta=0, \quad \cos \theta = \frac{E_1}{P_\varphi}.
\end{equation}
Since in this section we set $\vec{Q}\neq 0$, we obtain $P_\varphi \neq 0$ from \eqref{eq_Q12}.
Taking \eqref{eqth} into account, the equations defining the angle $\varphi$
can be represented as
\begin{equation}
\label{eqphi}
\frac{d \varphi}{d\tau}=\frac{P_\varphi}{E_2^2 + E_3^2}\left(E_2 \frac{\partial H}{\partial E_2} + E_3 \frac{\partial H}{\partial E_3} \right).
\end{equation}
The Poisson bracket of the coordinate $t$ with all variables of the reduced system is equal to
zero. Hence,
\begin{equation}
\label{eq_dt}
\frac{d t}{d\tau} = - \frac{\partial H}{\partial \mathcal{E}}.
\end{equation}
Thus, the problem reduces to an analysis of the reduced system \eqref{eq_EZ}, and the
evolution of the remaining coordinates is governed by \eqref{eqth}, \eqref{eqphi} and \eqref{eq_dt}.

 In the variables of the reduced system, inequality \eqref{eq_ner1} can be represented as
\begin{equation}
\label{eqUnon}
\begin{gathered}
u_\alpha u^\alpha=\frac{1}{a(r)}\left(\frac{\partial H}{\partial P_r} \right)^2 + r^2 \left( \frac{\partial H}{\partial E_2} \right)^2 + r^2 \left( \frac{\partial H}{\partial E_3} \right)^2 - a(r) \left(\frac{\partial H}{\partial \mathcal{E} } \right)^2<0.
\end{gathered}
\end{equation}

\begin{rem}
If one does not assume that relations \eqref{eq_Q12} are satisfied, then
to define the remaining cordinates one needs to find a cyclic coordinate corresponding to
the integral $F$.
It satisfies the condition
\begin{equation}
\label{eq_psi}
\{ \psi, F \}={\rm const},
\end{equation}
where the Poisson bracket is defined by \eqref{eq_J2}.
If $\sqrt{\ell}>P_\varphi$ and $\sqrt{\ell}>L_1$, one can define the angles
$$
\cos \alpha = \frac{P_\varphi}{\sqrt{\ell}}, \quad \cos \beta = \frac{L_1}{\sqrt{\ell}},
$$
where $\alpha, \beta \in(0, \pi)$. Then condition \eqref{eq_psi} is satisfied by $\psi$,
which is defined by the equation
\begin{equation}
\label{eq_cos}
\cos \theta = \cos \alpha \cos \beta + \sin \alpha \sin \beta \cos \psi.
\end{equation}
This equation has the following geometric interpretation. Define
a spherical triangle with sides $\alpha$, $\beta$, $\theta$. Then relation \eqref{eq_cos}
is a theorem of cosines where $\psi$ is the angle between sides $\alpha$ and $\beta$.
Previously, this cyclic coordinate $\psi$ and its geometric interpretation were
obtained in Ref. \cite{Ramond2025}.

The cyclic variable $\psi$ has the following nonzero Poisson brackets with
reduced variables:
$$
\{ \psi, E_2 \} = \frac{E_2}{\sqrt{\ell} \sin^2\beta}, \quad
\{ \psi, E_3 \} = \frac{E_3}{\sqrt{\ell} \sin^2\beta}.
$$
As a result, the cyclic coordinate satisfies the analogous equation \eqref{eqphi}.
In this approach, the angle $\varphi$ can be defined, for example, from the solution
of the equation $Q_1={\rm const}$.
 However, throughout the remainder of this paper we will assume for simplicity that
 relations \eqref{eq_Q12} hold.
\end{rem}

\subsection{Reduction in the case $\vec{Q} = 0$}
\label{sectionQ_zero}
Consider separately the case $\vec{Q} = 0$. It follows from \eqref{eq_Int} that
system \eqref{eq_HamSys} has the invariant manifold
\begin{equation}
\label{eq_PQ0}
P_\theta=0, \quad L_1=0, \quad P_\varphi=0.
\end{equation}
In this case, according to \eqref{eq_EZ0}, we have $\vec{E}=0$, i.e., for the reduced
system \eqref{eq_EZ} these trajectories correspond to the zero value of the Casimir
function $F(\vec{E})=0$.

Taking $Z_2>0$ into account, the Hamiltonian \eqref{eq_EZ} simplifies to
$$
\begin{gathered}
H_0 = \left. H \right|_{\vec{E}=0} = a(r)\left[ \frac{1}{2} + \frac{2 Z_1^2 + Z_2^2}{V_0} \right]P_r^2 + Z_2 Z_3\frac{Z_1 + 2 r \mathcal{E} }{r V_0}P_r +
\frac{Z_2^2}{r^3 V_0}\left[ r a(r) Z_2^2 + \mu Z_1^2\right] - \\
 \frac{(r \mathcal{E} + \mu Z_1)^2}{r^3(r - 2 \mu)}\left[ \frac{1}{2} + \frac{2 Z_1^2 - Z_3^2}{V_0} \right] + \frac{Z_2^2}{r} \left[ \mathcal{E} Z_1 + \frac{a(r)}{2r^2}
 - \frac{a(r)Z_3^2}{r V_0} \right], \\
  V_0=\frac{m^2}{\mu}r^3 - 2Z_1^2 - Z_2^2 + Z_3^2.
\end{gathered}
$$
Thus, the problem reduces to an analysis of a Hamiltonian system with two degrees of freedom
\begin{equation}
\label{eq_ZrPr}
\frac{d \vec{Z}}{d\tau} = ({\bf J}\vec{Z}) \times \frac{\partial H}{\partial \vec{Z}}, \quad
\frac{d P_r}{d\tau}= - \frac{\partial H}{\partial r} ,\quad
\frac{d r}{d\tau} = \frac{\partial H}{\partial P_r}.
\end{equation}
The trajectories of this system lie on the symplectic leaf
$$
\mathcal{S}_c^4= \{ (\vec{Z}, r, P_r) \ | \ C_\circ(\vec{Z})=c^2\}.
$$
The Tulczyjew condition defines the invariant manifold
\begin{equation}
\label{tyleq}
 Z_1 P_r + \frac{Z_2Z_3}{r}=0, \quad \mathcal{E} Z_1 + \frac{\mu Z_1^2}{r^2} - \frac{a(r)}{r}Z_2^2=0.
 \end{equation}
On this invariant manifold, the system \eqref{eq_ZrPr} is integrable by quadratures.
Indeed, from Eqs. \eqref{tyleq}, the Casimir function $C_\circ(\vec{Z})=c^2$
and relation \eqref{eq_mmm} (taking $\vec{E}=0$ into account), we obtain two solutions
$$
\begin{gathered}
\mathcal{R}_{\pm} = \left\{ Z_1= \frac{\mathcal{E} c^2 r^2}{m^2r^3 - \mu c^2}, \ Z_2 = \frac{r^3c\sqrt{a(r)}m\mathcal{E} }{m^2 r^3 - \mu c^2}, \
Z_3=\pm \frac{c\sqrt{N(r)}}{\sqrt{a(r)} (m^2 r^3 - \mu c^2)}, \ P_r = \pm\frac{m \sqrt{N(r)}}{a(r)(m^2 r^3 - \mu c^2)} \right\}, \\
N(r)=r^4(m^2r^2 - c^2a(r))\mathcal{E}^2 - a(r)(m^2r^3 - \mu c^2)^2,
\end{gathered}
$$
where the upper sign corresponds to $\mathcal{R}_+$, and the lower sign, to $\mathcal{R}_-$.
Note that the choice of
variables \eqref{eq_EZ0} leads to $Z_2>0$.
Hence, for the solutions $\mathcal{R}_{\pm}$ it is assumed that the condition $Z_2>0$ is
always satisfied (otherwise one should choose two other solutions which differ in that the
sign reversal $Z_2 \to -Z_2$ occurs).

Substituting the resulting solutions into the last of equations \eqref{eq_ZrPr}, we obtain
the following quadrature for defining the evolution of the radial coordinate:
\begin{equation}
\label{eqdrr}
\begin{gathered}
\frac{dr}{d\tau} = \pm \frac{\sqrt{N(r)}}{\eta(r)}, \quad
\eta(r) = m^2 r^3 - \mu c^2 - \frac{3\mu c^4r^4 \mathcal{E}^2}{(m^2 r^3 - \mu c^2)^2},
\end{gathered}
\end{equation}
where the upper sign corresponds to $\mathcal{R}_+$, and the lower sign, to $\mathcal{R}_-$.

Now consider how one can retrieve the values of the coordinates $\theta$, $\varphi$ and $t$
from the known solution of \eqref{eq_ZrPr}.
The nonzero Poisson brackets \eqref{eq_J2} of these coordinates with the variables of the
reduced system have the form
$$
\begin{gathered}
\{\theta, E_2 \}=\frac{L_3}{\sqrt{L_2^2 + L_3^2}}, \quad \{\theta, E_3 \}=-\frac{L_2}{\sqrt{L_2^2 + L_3^2}}, \\
\{\varphi, E_2 \}=-\frac{L_2}{\sin \theta \sqrt{L_2^2 + L_3^2}}, \quad \{\varphi, E_3 \}=-\frac{L_3}{\sin \theta \sqrt{L_2^2 + L_3^2}}.
\end{gathered}
$$
It turns out that system \eqref{eq_HamSys} on the invariant manifold described by relations
\eqref{eq_PQ0} and by the value of the Casimir function $C_\star=0$ has the additional
integral
$$
K=\frac{L_2 \sin \theta}{\sqrt{L_2^2 + L_3^2}}.
$$
Let us fix the value of this integral $K(\vec{L}, \theta)=\kappa<1$ and express $L_2$ and
$L_3$ in terms of $\kappa$ and $Z_2$ as follows:
\begin{equation}
\label{eq_L23}
L_2=\frac{\kappa Z_2}{\sin \theta}, \quad L_3=\pm Z_2 \sqrt{1 - \frac{\kappa^2}{\sin^2 \theta}}.
\end{equation}
With this in mind, the equations can be written as
$$
\begin{gathered}
\frac{d\theta}{d\tau} = \pm \sqrt{1 - \frac{\kappa^2}{\sin^2 \theta}} \left. \frac{\partial H}{\partial E_2} \right|_{\vec{E}=0} - \frac{\kappa}{ \sin \theta} \left. \frac{\partial H}{\partial E_3} \right|_{\vec{E}=0}, \\
\frac{d\varphi}{d\tau} = - \frac{\kappa}{ \sin^2 \theta} \left. \frac{\partial H}{\partial E_2} \right|_{\vec{E}=0} \mp \frac{1}{\sin \theta}\sqrt{1 - \frac{\kappa^2}{\sin^2 \theta}} \left. \frac{\partial H}{\partial E_3} \right|_{\vec{E}=0}, \\
\frac{dt}{d\tau} = - \frac{\partial H_0}{\partial \mathcal{E}},
\end{gathered}
$$
where the upper sign corresponds to the case $L_3>0$, and the lower sign, to the case
$L_3<0$. Each time the radicand in \eqref{eq_L23} vanishes, the upper and lower signs in
these equations must be interchanged.

Finally, taking $\left. \frac{\partial H}{\partial E_3} \right|_{\vec{E}=0}=0$
and the relations $\mathcal{R}_{\pm}$ into account, we obtain
\begin{equation}
\label{eqQ_0_th_phi}
\begin{gathered}
\frac{d\theta}{d\tau} = \pm\frac{\Phi(r)}{\eta(r)}\sqrt{1 - \frac{\kappa^2}{\sin^2 \theta}}, \quad \frac{d\varphi}{d\tau} = - \frac{\kappa \Phi(r)}{\eta(r) \sin^2 \theta},
\quad \Phi(r) = - c \mathcal{E} r \frac{m^2r^3 + 2\mu c^2}{m^2 r^3 - \mu c^2}, \\
\frac{dt}{d\tau} = \frac{m r^3 \mathcal{E} }{a(r)}.
\end{gathered}
\end{equation}

It is important to emphasize that the case $Z_2=0$ requires a separate analysis since
$\vec{L}=0$ and, according to \eqref{eq_EZ0}, the coordinates $Z_1$ and
$Z_2$ are undefined. However, in the limit $c \to 0$ the resulting equations \eqref{eqQ_0_th_phi} and \eqref{eqdrr} 
become equations of the radial geodesics of the
Schwarzschild metric \cite{chandrasekhar1998mathematical}. Indeed, in the limit $c \to~0$ we
have $\Phi(r)=0$, whence $\dfrac{d\theta}{d\tau} = \dfrac{d\varphi}{d\tau} =0$ and
 the equation for the radial coordinate reduces to the well-known equation \cite{chandrasekhar1998mathematical}
 $$
\left. \left( \frac{dr}{d\tau} \right)^2 \right|_{c \to 0} = \frac{2\mu}{r} - 1 + \frac{ \mathcal{E}^2}{m^2}.
 $$
If $c \neq 0$, the world line is, according to \eqref{eqQ_0_th_phi}, not a radial trajectory.
 A qualitative analysis of the trajectories of system \eqref{eq_ZrPr} is a separate
 problem and is therefore not carried out here. In what follows, we will consider only
the case $\vec{Q}\neq 0$.

\begin{rem}
Refs. \cite{Mukherjee2022, Costa2015} describe trajectories for which
the rotating test body falls radially
onto the event horizon, following the geodesic. Such trajectories are fixed points of
the symmetry field \eqref{eq_symu} for which $P_\theta=0$, $L_2=L_3=0$, $M_2=M_3=0$ and
$L_1 \cos \theta - P_\varphi =0$. In this situation, satisfying the Tulczyjew condition
implies that $M_1 = 0$. As we see, for these fixed points,
the reduced coordinates \eqref{eq_EZ0} turn out to be undefined, and so these trajectories
are not described by the reduced systems obtained in Sections \ref{sectionQ_zero},
\ref{sectionRedQneq0} and are not considered here. Note that, for these trajectories, $F = C_\circ \neq 0$. 
\end{rem}

\subsection{Canonical variables}
\label{eqCanonic}
For the resulting reduced system \eqref{eq_EZ}, the Poisson bracket \eqref{SPa} is not canonical.
The canonical variables introduced for this bracket are closely related to the variables
used in Ref. \cite{Ramond2025}.
We consider this relationship in more detail.

With the notation used here, the variables $(\zeta, \pi_\zeta, s, \pi_s)$ from
Ref. \cite{Ramond2025} are given as follows:
\begin{equation}
\label{eqRamond1}
\begin{gathered}
L_1 = \pi_s, \quad L_2 = \sqrt{\pi_\zeta^2 - \pi_s^2} \cos \sigma, \quad L_3 = \sqrt{\pi^2_\zeta - \pi_s^2}\sin \sigma, \\
M_1 = - \sqrt{1 - \frac{c^2}{\pi_\zeta^2}} \sqrt{\pi_\zeta^2 - \pi_s^2}\sin \zeta, \quad
M_2 = \sqrt{1 - \frac{c^2}{\pi_\zeta^2}}(\pi_s \sin \zeta \cos \sigma + \pi_\zeta \cos\zeta \sin \sigma), \\
M_3 = \sqrt{1 - \frac{c^2}{\pi_\zeta^2}}(\pi_s \sin \zeta \sin \sigma - \pi_\zeta \cos\zeta \cos \sigma),
\end{gathered}
\end{equation}
where the angle $\sigma$ is defined from the equation
$$
\tan(\sigma - s) = \frac{P_\theta \sin \theta}{ \pi_s \cos \theta - P_\varphi }.
$$
These variables explicitly parameterize the symplectic leaf $\mathcal{M}^{12}_c$.
Substituting \eqref{eqRamond1} into the definition of the reduced variables \eqref{eq_EZ0}
and then transforming the resulting relations using \eqref{eq_ES6c}, we find
$$
\begin{gathered}
E_1 = \pi_s, \quad E_2 = -\sqrt{\ell - \pi_s^2} \cos s, \quad E_3 = -\sqrt{\ell - \pi_s^2} \sin s,\\
Z_1 = \sqrt{\pi_\zeta^2 - c^2} \sin \zeta, \quad Z_2=\pi_\zeta, \quad Z_3 = \sqrt{\pi_\zeta^2 - c^2}\cos \zeta.
\end{gathered}
$$
Straightforward calculations show that the Poisson bracket \eqref{SPa} in the new variables is
canonical
since the nonzero Poisson brackets have the form
$$
\{ \zeta, \pi_\zeta \} = 1, \quad \{ s, \pi_s \}=1.
$$

Ref. \cite{Ramond2025} is focused on obtaining a reduced Hamiltonian
system linearized in the spin, and is an extension of Ref. \cite{Ramond2022}. A special
feature of this system is that it has the additional integral
$r M_1 = {\rm const }$ and hence is integrable by the Liouville\,--\,Arnold theorem.
The solution in the form of elliptic Jacobi functions for this problem was obtained in
Ref. \cite{Witzany2024} using the Marck tetrad formalism.
Separation of variables to linear order in the spin for a Kerr metric is performed in
Ref. \cite{Skoup2025}.

Refs. \cite{Ramond2025, Ramond2022} do not derive an explicit form of the Hamiltonian
that is not linearized in the spin. Therefore, the reduced equations of motion obtained in
Section \ref{sectionRedQneq0} or those equivalent to them are absent in Refs.
\cite{Ramond2025, Ramond2022}. Note that the nonlinear system \eqref{eq_Hi}-\eqref{eq_EZ}, as opposed to the system linearized in the spin is not integrable by quadratures. Therefore, the
relationship between two systems is a separate problem (see the section ``Conclusions and
Prospects'' in Ref. \cite{Ramond2022}), and so it will not be discussed further.

In Ref. \cite{Witzany2019}, a different set of canonical variables is introduced.
It also parameterizes explicitly the symplectic leaf $\mathcal{M}^{12}_c$
and is defined by the relations
\begin{equation}
\label{eqWitzany2019}
\begin{gathered}
\varphi_W=\arctan \frac{L_2}{L_1}, \quad \psi_W = \arctan \frac{L_2}{L_1} - \kappa_W \arccos\big( M_2\sqrt{\mathscr{C}} \big), \\
A_W = L_3 - \sqrt{(\vec{L}, \vec{L})}, \quad B_W = \sqrt{(\vec{L}, \vec{L})} - c, \\
\kappa_W={\rm sign}(M_2L_1 - M_1L_2), \quad \mathscr{C} = \frac{(\vec{L}, \vec{L})}{(L_1^2 + L_2^2)(\vec{M}, \vec{M})}.
\end{gathered}
\end{equation}
For the variables $(\varphi_W, A_W, \psi_W, B_W)$, the Poisson structure \eqref{eq_J2}
reduces to a canonical one since the nonzero Poisson brackets have the form
$$
\{ \varphi_W, A_W \} = 1, \quad \{ \psi_W, B_W \}=1.
$$
The variables $\varphi_W$, $\psi_W$ and $B_W$ are not integrals of the symmetry field
\eqref{eq_symu}. Therefore, after changing to these variables in the system
\eqref{eq_HamSys}, the Hamiltonian $H$ depends explicitly on $\theta$ and $P_\theta$,
in contrast to the variables $\vec{E}$ and $\vec{Z}$ (see equation \eqref{eq_Hi}).

In other words, the variables \eqref{eqWitzany2019} do not define reduction by the 
symmetry group $SO(3)$ (for the concept of reduction see, e.g., the review 
\cite{Borisov2015Symmetries}). Therefore, in these variables the
system \eqref{eq_HamSys} for the Schwarzschild metric does not reduce to analysis of a 
two-dimensional Poincar\'{e} map, in contrast to the system \eqref{eq_EZ}(see Section \ref{sectionMap}). 
In order to reduce the system to analysis of a two-dimensional map, 
the authors of \cite{Witzany2019} suggested restricting the system to an invariant 
manifold that defines special planar motions.
Note that Eqs. \eqref{eqth} can be used as an invariant manifold, but it should be kept in mind 
that Eqs. \eqref{eqth} do not define the Poisson submanifold for the bracket
\eqref{eq_J2}, and so the Hamiltonian of the system and the Poisson bracket cannot be
restricted to these relations. To do so, it is necessary to project the Poisson bracket
onto the submanifold defined by relations \eqref{eq_Q12}, which leads to the Dirac bracket \cite{hanson1976constrained, ArnoldKozlov}.

Of course, the application of canonical variables has a number of advantages because 
it is for these variables that the methods of Hamiltonian 
mechanics (KAM theory etc.) were developed.
However, throughout the remainder of this paper, the reduced variables 
$\vec{E}$ and $\vec{Z}$ will be used.  They turn out to be more convenient to search for 
relative equilibria because in these redundant coordinates
the reduced equations of motion preserve the algebraic form (i.e., they do not contain 
various trigonometric functions).

\section{Relative equilibria}
\label{section_eqP}
We describe the fixed points of the reduced system \eqref{eq_EZ} obtained above,
which are relative equilibria. For them $r={\rm const}$ and $E_1={\rm const}$. Then
it follows from \eqref{eqth} that $\theta ={\rm const}$. Therefore, for relative equilibria
the test body moves generally in a circular orbit.

In the reduced system \eqref{eq_EZ}, one can assume without loss of generality that
\begin{equation}
\label{eq_m_mu}
m=1, \quad \mu=1.
\end{equation}
Indeed, the mass of the test body and the black hole can be excluded from system
\eqref{eq_EZ} after transforming from the initial variables and energy to the
dimensionless ones
$$
P_r \to m \overline{P}_r, \ r \to \mu \overline{r}, \ E_i \to m \mu \overline{E}_i, \ Z_i \to m \mu \overline{Z}_i, \
\mathcal{E} \to m \overline{\mathcal{E}}.
$$
After this change of variables, the dimensionless total angular momentum has the form
$$
\overline{c}=\frac{c}{m\mu}.
$$
The most realistic test body is one of the following compact objects \cite{Hartl2003}:
\begin{itemize}
\item [---] a black hole that is described by a Kerr metric and has a total angular momentum
that satisfies $c< m^2$;
\item [---] a neutron star for which, according to \cite{Cook1994}, various state equations
predict $c\leqslant 0.6 m^2$.
\end{itemize}
Taking into account the assumption about the smallness of the mass of the test
body $m \ll \mu$, we obtain $\overline{c}~\ll~1$.
However, we will consider the case $\overline{c} \sim 1$ to completely analyze the
dynamics of the MPD equations and to supplement the results of \cite{Suzuki1998, Hackmann2014}.

In addition, formally one can obtain $\overline{c} \sim 1$
for naked singularities, which arise for some solutions of the Einstein equations.
Naked singularities have not been observed in nature, and their existence would violate
the Cosmic Censorship Hypothesis, but they continue to be the subject of intense research
(see, e.g., \cite{Nguyen2023}).

Next, we will assume that conditions \eqref{eq_m_mu} are satisfied since this
will allow us to simplify the explicit relations presented below.

\subsection{Symmetric circular orbits}
For the reduced system \eqref{eq_EZ} the symmetry \eqref{eq_Sym} is preserved and has the form
$$
\mathcal{I}^{(r)}: \ E_1\to -E_1, \quad E_3\to -E_3,
$$
where the other variables remain unchanged.

We fix the coordinates on the line of the fixed points of this symmetry
$$
 {\rm Fix } \ \mathcal{I}^{(r)}: \ E_1=0, \ E_3=0.
 $$
Then from the invariant submanifold \eqref{eq_Tred},
the symplectic leaf $\mathcal{S}^6_c$ and the integral \eqref{eq_mmm} we obtain two solutions:
$$
\begin{gathered}
\mathcal{N}_0^\pm = \left\{ (\vec{E}, \vec{Z}, r, P_r) \ | \ \vec{E}=(0, \pm\sqrt{\ell}, 0), \ \vec{Z}=\vec{Z}^{\pm}, \ P_r=P_r^\pm \right\}. \\
\end{gathered}
$$
Here and in what follows, the upper sign corresponds to $\mathcal{N}_0^+$, and the lower sign,
to $\mathcal{N}_0^-$. The vectors $\vec{Z}^{+}$, $\vec{Z}^{-}$ and the linear momentum
$P_r^\pm$ are given by the following relations:
$$
\begin{gathered}
Z^{\pm}_1=-\frac{r^2 c ( \mathcal{E} c \mp \sqrt{\ell} )}{ r^3 - c^2 }, \quad Z^{\pm}_2=\frac{c\sqrt{b(r)}(r^3 \mathcal{E} \mp c \sqrt{\ell} ) }{ r^3 - c^2}, \\
Z^{\pm}_3 = c\sqrt{b(r)} P_r^\pm, \\
\frac{( r^3 \mathcal{E} \mp c \sqrt{\ell} )^2}{b(r)(r^3 - c^2)^2} - \frac{r^4 ( \mathcal{E} c \mp \sqrt{\ell})^2}{(r^3 - c^2)^2} - b(r)( P_r^{\pm})^2 = 1,
\end{gathered}
$$
where the following notation has been introduced:
$$
b(r)=1 - \dfrac{2}{r}.
$$
After substituting these values of the reduced variables into the equation for the radial
coordinate, we obtain
\begin{equation}
\label{eq_Rpm}
\begin{gathered}
\left(\frac{dr}{d\tau} \right)^2 = \frac{R_\pm(r)}{\chi_\pm^2(r)}, \\
R_\pm(r) = ( \mathcal{E} r^3 \mp c\sqrt{\ell} )^2 - r^4 b(r) (\mathcal{E} c \mp \sqrt{\ell})^2 - b(r)( r^3 - c^2)^2, \\
\chi_\pm(r) = r^3 - c^2 - 3 \left[ \frac{r^2c (\mathcal{E}c \mp \sqrt{\ell}) }{r^3 - c^2 } \right]^2.
\end{gathered}
\end{equation}
Previously, these relations were obtained, using different notation, independently in
Refs. \cite{Tod1976, Saijo1998}. Moreover, these papers present derivations of
equations on the invariant manifold \eqref{eq_th_p_S} for
a more general case of the Kerr metric.

We represent Eqs. \eqref{eqphi} and \eqref{eq_dt} on $\mathcal{N}_0^\pm$ as
$$
\begin{gathered}
\frac{d\varphi}{d\tau} = \frac{r(r^3 + 2c^2)(c \mathcal{E} \pm \sqrt{\ell} )}{(r^3-c^2)\chi_\pm(r)}, \quad
\frac{dt}{d\tau} = \frac{r(r^3\mathcal{E} \mp c \sqrt{\ell} )}{(r - 2)\chi_\pm(r)}.
\end{gathered}
$$
Inequality \eqref{eqUnon} reduces to
\begin{equation}
\label{eq_Phi_U}
u_\alpha u^\alpha = \frac{(r^3 + 2c^2)(r^3 - c^2) -(2r^3 + c^2)\chi_\pm(r)}{ \chi^2_\pm(r)}<0.
\end{equation}
The behavior of the radial coordinate in Eq. \eqref{eq_Rpm} is examined in detail in
Ref. \cite{Hackmann2014}. We consider only the case of circular orbits $r={\rm const}$,
in which $\mathcal{N}_0^\pm$
define the fixed points of the reduced system \eqref{eq_EZ}. In this case, $P_r=0$ and $r$ is
the critical point of the function $R_\pm(r)$:
\begin{equation}
\label{eq_R_critic}
R_\pm(r)=0, \quad \frac{dR_\pm(r)}{dr} =0.
\end{equation}
For such circular orbits the angular momentum of the body is parallel to
the total angular momentum. Note that for $\mathcal{N}^+_0$ they are aligned,
and for $\mathcal{N}^-_0$ they have opposite directions.
\begin{figure*}[!ht]
\begin{center}
\includegraphics[scale=0.9]{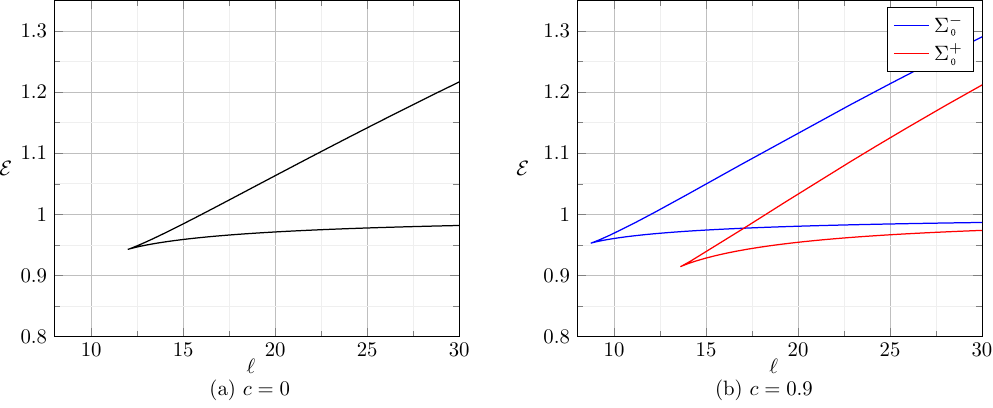}
\end{center}
\caption{ Curves $\Sigma_0^\pm$ on the plane $(\ell, \mathcal{E})$ for different $c$.}
		\label{fig4}
\end{figure*}

\begin{figure*}[!ht]
\begin{center}
\includegraphics[scale=0.9]{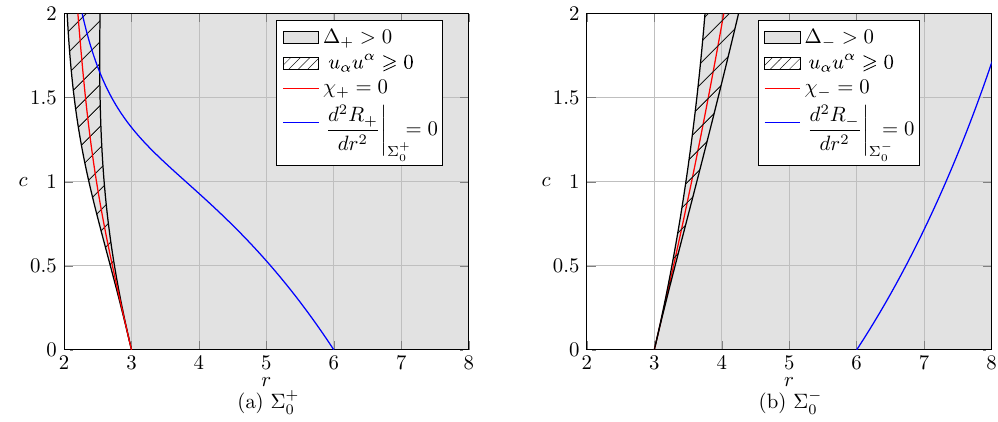}
\end{center}
\caption{Domain of definition (gray) of the radial coordinate $r$  for symmetric 
circular orbits depending on
$c$ for $\Sigma^{\pm}_0$.  Light blue indicates the curve separating 
circular orbits that are stable or unstable  with respect to perturbations lying on the invariant 
manifold
$ {\rm Fix } \ \mathcal{I}^{(r)}$.} 
\label{fig5}
\end{figure*}

The system of equations \eqref{eq_R_critic} has a rather complex dependence on the radial
coordinate $r$, whereas the terms $\mathcal{E}$ and $\sqrt{\ell}$ appear in Eqs. \eqref{eq_R_critic}
only in the first and the second degree.
Solving this system for the above values of the first integrals, we obtain
\begin{equation}
\label{eq_Sigma_pm}
\begin{gathered}
\Sigma^{\pm}_{0} = \left\{ (\ell, c, \mathcal{E}) \ | \ \ell=\frac{2 \eta^2}{r \Delta_\pm}, \ \mathcal{E}=\mathcal{E}_\pm(r, c) \right\}, \\ \eta = r^6 -2c^4 (r - 2)^2 - c^2r^4 (r - 7) - 11 c^2r^3\\
\mathcal{E}_\pm = \pm\frac{ 2\big((2r - 3)c^4 - 3c^2r^3 + r^6(r - 3)\big)\sqrt{\dfrac{\eta^2}{\Delta_\pm}} }{r^4c\sqrt{2r}(2r^3 + c^2 - 9r^2)}
\mp \frac{ (r^3 - c^2)\sqrt{\Delta_\pm}}{r^4c\sqrt{2r}(2r^3 + c^2 - 9r^2)},
\end{gathered}
\end{equation}
where $\Sigma_0^+$ corresponds to $\mathcal{N}_0^+$, $\Sigma_0^-$ corresponds to
$\mathcal{N}^-_0$, and the following function has been introduced:
$$
\begin{gathered}
\Delta_\pm = 2r^3(r - 3)(r^3 - c^2)^2 \mp r^2c(r - 2)(2r^3 + c^2 - 9r^2) \sqrt{4r^7 + 13c^2r^4 - 8c^4r} + \\ c^2(r - 2)\big[(r - 2)c^2\big(4(2r - 3)c^2 + r^3(8r - 39)\big) + 2r^7(r - 8) + 33r^6 \big].
\end{gathered}
$$
To describe the surfaces $\Sigma^+_0$ and $\Sigma^-_0$, it is convenient to consider
the section formed by their intersection with the plane $c={\rm const}$ on which they
define two curves.

If $c=0$, then we obtain the case of geodesics in which relations \eqref{eq_Sigma_pm}
define one curve (see Fig. \ref{fig4}a) on which the radial coordinate changes in the interval
$r\in(3, +\infty)$. The value $r=6$ corresponds to a singular cusp point on this curve
(see, e.g., \cite{Stewart2006}), which separates stable
circular orbits with $r\in(6, +\infty )$ from unstable orbits.

Further, as the total angular momentum $c$ increases, the curves defined by
\eqref{eq_Sigma_pm} no longer coincide (see Fig. \ref{fig4}b), but the singular cusp point
persists \cite{Jefremov2015, Hackmann2014}. Also, on each curve the domain of definition
of the radial coordinate $r$ is given by the condition
$\Delta_\pm>0$ for $\Sigma^{\pm}_{0}$, respectively. The value of the radial coordinate
$r$ defining the cusp satisfies the equation
 $$
 \left. \frac{d^2R_\pm(r)}{dr^2} \right|_{\Sigma^{\pm}_{0}} =0.
 $$
How the domain of definition of the radial coordinate and the position of the
cusp change depending on $c$ is illustrated in Fig. \ref{fig5} separately for each of the
curves, $\Sigma_0^{+}$ and $\Sigma_0^{-}$.
The region in which inequality \eqref{eq_Phi_U} is not satisfied  (i.e., 
the velocity ceases to be a time-like vector) is hatched; hence, this region contains 
nonphysical circular orbits.
As we see, in this domain the functions $\chi_\pm(r)=2r^{-3}\left.(d + 2p_\alpha p^\alpha)\right|_{\Sigma_0^\pm}$
can vanish (see Remark \ref{rem2}), but the pole-dipole approximation cannot be applied any longer, and so the
corresponding trajectories should be analyzed by taking higher multipole moments into account.

It is also seen from Fig. \ref{fig5} that, due to the spin-orbit interaction, for
$\mathcal{N}_0^+$ the circular orbits can be nearer to the event horizon. Conversely,
for $\mathcal{N}_0^-$ they can be farther from the event horizon (for details, see
\cite{Suzuki1997, Jefremov2015}).

In the general case, for rotating test bodies on the curve $\Sigma_0^{+}$
the singular cusp point separates trajectories that are stable and unstable only
with respect to perturbations lying on the invariant manifold
$ {\rm Fix } \ \mathcal{I}^{(r)}$. It was shown in Ref. \cite{Suzuki1998} that
the equilibrium points considered can lose stability with respect to perturbations that
do not lie in $ {\rm Fix } \ \mathcal{I}^{(r)}$. It turns out that this is closely related
to the birth of other asymmetric relative equilibria (i.e., those that do not lie in $ {\rm Fix } \ \mathcal{I}^{(r)}$). Previously, these equilibrium points were not found, and so we consider
them separately.

\subsection{Asymmetric circular orbits}
\label{nonsymmetricOrbit}
On the invariant submanifold \eqref{eq_Tred} and at a fixed value of the integral
\eqref{eq_mmm}, the reduced system \eqref{eq_EZ} possesses
two one-parameter families of equilibrium points
\begin{equation}
\label{eq_N1}
\begin{gathered}
\mathcal{N}_1^\pm = \left\{ (\vec{E}, \vec{Z}, r, P_r) \ | \ \vec{E}=(0, E^{(1)}_2, E_3^{\pm} ), \right.\\ \left. \vec{Z}=(z, \sqrt{z^2 + 27(1 - z \mathcal{E})}, 0), \ r=3, \ P_r=0 \right\}, \\
E_2^{(1)}=\frac{12\sqrt{3}(3 - 2 z \mathcal{E})) }{\sqrt{z^2 + 27(1 - z \mathcal{E})}}, \ E_3^{\pm}=\pm3\sqrt{\frac{(1- z \mathcal{E})(z - 9 \mathcal{E})^2 }{z^2 + 27(1 - z \mathcal{E})} - 1},
\end{gathered}
\end{equation}
where $z$ is a parameter. As we see, these equilibrium points are, according to \eqref{eqth},
circular orbits lying in the equatorial plane $\theta=\dfrac{\pi}{2}$ for which
$$
\frac{d\varphi}{d\tau} = \frac{P_\varphi}{9}, \quad \frac{dt}{d\tau} = 3 \mathcal{E}.
$$

\begin{figure*}[!ht]
\begin{center}
\includegraphics[scale=0.8]{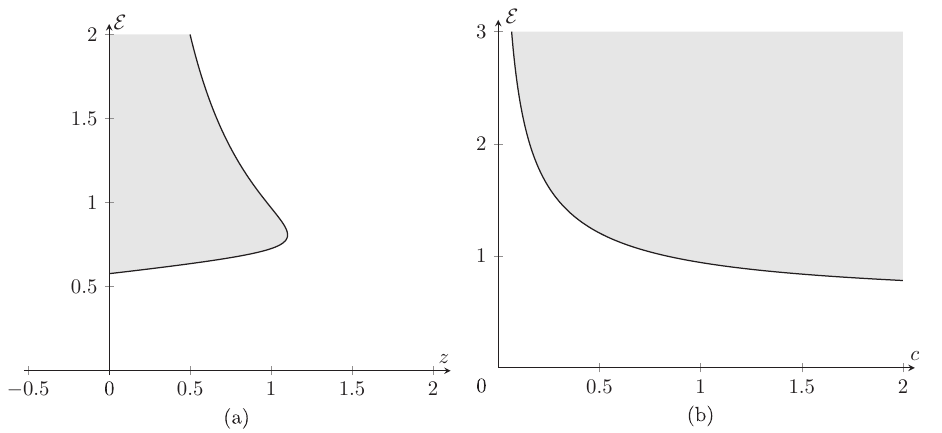}
\end{center}
\caption{Region on the plane $(z, \mathcal{E})$ and $(c, \mathcal{E})$ for which
there exist equilibrium points $\mathcal{N}_1^\pm$ (gray). Black indicates the curve on
which $E_3^\pm=0$.}
\label{fig6}
\end{figure*}

To describe the domain of definition of the parameter $z$, we note that for
$\mathcal{N}_1^\pm$
inequality \eqref{eqUnon} is equivalent to
$$
u_\alpha u^\alpha =-z \mathcal{E}<0.
$$
The domain of definition of $z$ is given by this inequality and by the condition for
positivity of the expressions under
the square root in $\mathcal{N}_1^\pm$.
This region on the plane $(z, \mathcal{E})$ is shown in Fig. \ref{fig6}a.

We express $z$ in terms of the value of the Casimir function $C_\circ(\vec{Z})=c^2$:
$$
z= \frac{27 - c^2}{27\mathcal{E}},
$$
the corresponding region is shown in Fig. \ref{fig6}b.
Then, taking this relation into account and substituting \eqref{eq_N1} into $F(\vec{E})=\ell$,
we obtain
$$
\Sigma_{1} = \left\{ (\ell, c, \mathcal{E}) \ | \ \frac{c^2}{3} + 27\mathcal{E}^2 - \ell=9 \right\},
$$
where $\ell \in(\ell^*, +\infty)$ and $$
\ell^* = \frac{2}{9}c^2 - \frac{3}{2} + \frac{3}{c}\sqrt{27 + \frac{13}{4}c^2 - \frac{2}{27}c^4}.
$$

The equilibrium points $\mathcal{N}_1^\pm$ have $E_1=0$. Therefore,
they lie, according to \eqref{eqth}, in the equatorial plane $\theta = \dfrac{\pi}{2}$.
In addition to them, the reduced system \eqref{eq_EZ} also possesses equilibrium points
that have $E_3=0$, but $E_1 \neq 0$.
According to \eqref{eqth}, for them it holds that $\theta \neq \dfrac{\pi}{2}$ and they are not given by explicit relations, and so
we represent them in the following implicit form:
\begin{equation}
\label{eq_N2}
\begin{gathered}
\mathcal{N}_2^\pm = \left\{ (\vec{E}, \vec{Z}, r, P_r) \ | \ \vec{E}=(\pm\sqrt{z^2 - W^2}, E^{(2)}_2, 0 ), \ \vec{Z}=(\sqrt{z^2 - c^2}, z, 0), \ P_r=0 \right\}, \\
E_2^{(2)}=W \left( \sqrt{b(r)} - \frac{z^2 - c^2}{rz^2\sqrt{b(r)}} \right) + \frac{r\mathcal{E}\sqrt{z^2 - c^2}}{z\sqrt{b(r)}}, \ W=\frac{r^2 z}{c\sqrt{z^2 - c^2}}\big(
c \mathcal{E}\ - z \sqrt{b(r)}\big),
\end{gathered}
\end{equation}
where the parameters $r$ and $z$ satisfy the following system:
\begin{equation}
\label{eq001}
\begin{gathered}
6 r^\frac{7}{2}c^3\mathcal{E} - \sqrt{3 \zeta(z^2 - c^2)} + 3 z\sqrt{r-2} \big(2c^2 r^3 -(2r^3 + c^2)(z^2 - c^2) \big)=0,\\
\frac{z\sqrt{3 \zeta }}{2} \left[ z^2(r - 2)\big(c^4 - (2r^3 + c^2)z^2\big) + c^2 r^3 \left((r - 3) z^2 + c^2\left(r - \frac{5}{3}\right)\right) \right] + \\
\sqrt{(r-2)(z^2 - c^2)}\left( \frac{z^2\zeta}{2} + \frac{c^2 r^3 z^2}{2} (2r^3 + c^2) (3(z^2 - c^2)(r - 1) - 2c^2) + c^4r^6(c^2 + 3z^2)\right)=0, \\
\zeta=3z^2(r-2)(2r^3 + c^2)^2(z^2 - c^2) + 4c^4r^3(r^3 + 2c^2)>0.
\end{gathered}
\end{equation}

\begin{rem}
To obtain relations \eqref{eq001}, one needs to restrict system \eqref{eq_EZ} to
relations \eqref{eq_N2}.
This leaves us with the following two equations:
\begin{equation}
\label{eq002}
\left.\left( E_2\frac{\partial H}{\partial E_1} - E_1\frac{\partial H}{\partial E_2} \right) \right|_{\mathcal{N}_2^\pm}=0, \quad
\left. \frac{\partial H}{\partial r} \right|_{\mathcal{N}_2^\pm} =0.
\end{equation}
Next, by choosing a linear combination of these equations one can obtain a quadratic equation
in $\mathcal{E}$. Choosing the principal root of this equation, one obtains the first of
relations \eqref{eq001}. The other root must be excluded since inequality \eqref{eqUnon}
is not satisfied for it. Next, substituting the resulting solution for
$\mathcal{E}$ into any of Eqs. \eqref{eq002}, one can show that it reduces to the second of
Eqs. \eqref{eq001}. The computations themselves are rather cumbersome and so they are
not presented here explicitly.
\end{rem}

\begin{figure*}[!ht]
\begin{center}
\includegraphics[scale=0.65]{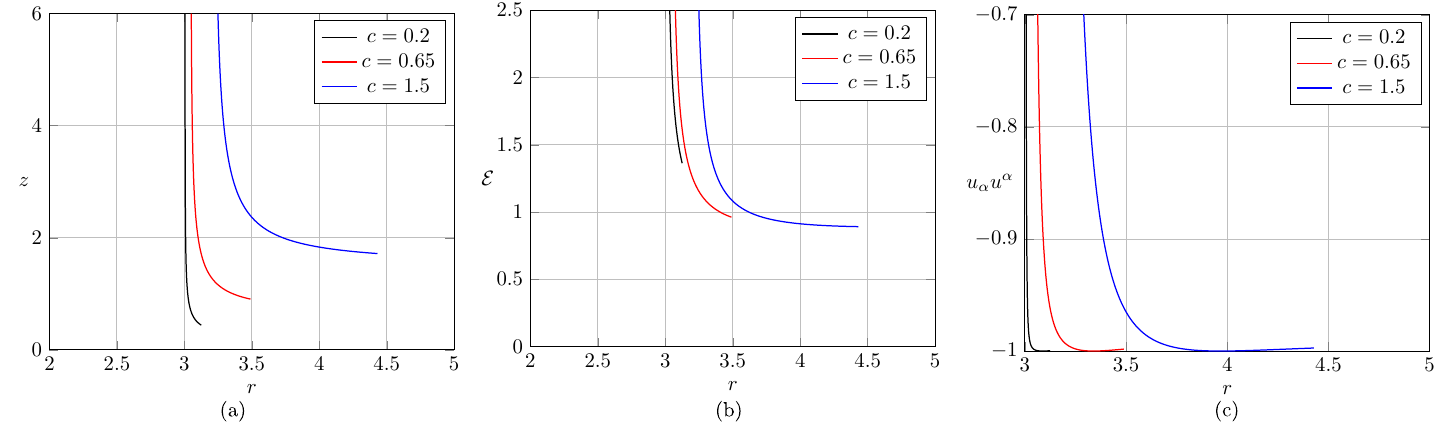}
\end{center}
\caption{Dependence $z(r)$, $\mathcal{E}(r)$ and $u_\alpha u^\alpha(r)$ for the equilibrium points $\mathcal{N}_2^\pm $
for different values of $c$.}
\label{fig7}
\end{figure*}

Figure \ref{fig7} shows how the dependences $z(r)$ and $\mathcal{E}(r)$ change for
the equilibrium points $\mathcal{N}_2^\pm $.
In these figures, the curves have an origin in which $z=W$. This corresponds to $E_1=0$.
As can be seen from Fig. \ref{fig7}c, for the equilibrium points 
$\mathcal{N}_2^\pm$, the condition $u_\alpha u^\alpha<0$ is satisfied, and for the 
equilibrium points $\mathcal{N}_2^\pm $, we have $r>3$.
Substituting these values of the coordinates for \eqref{eq_N2} into $F(\vec{E})=\ell$ and
taking \eqref{eq001} into account, we can write in implicit form
$$
\Sigma_{2} = \left\{ (\ell, c, \mathcal{E}) \ | \ \ell=\widetilde{\ell}(c, \mathcal{E}) \right\}.
$$

According to Figs.\ref{fig6}b and \ref{fig7}b, the circular orbits \eqref{eq_N1} and
\eqref{eq_N2} have energy $\mathcal{E}>1$ for small $c$. In the limit $c\to 0$ we have
$\mathcal{E} \to +\infty$, i.e., the asymmetric circular orbits found above arise only
when the angular momentum of the test body is nonzero.
We analyze how the angular momentum is directed in this case using the Pauli\,--\,Lubanski
vector:
$$
S^\alpha = \frac{\varepsilon^{\alpha \beta \mu \nu}}{2m \sqrt{-g}}p_\beta S_{\mu \nu},
$$
which satisfies the condition $p_\alpha S^\alpha=0$. In the reduced variables the
components of this vector can be represented as
$$
\begin{gathered}
S^t=\frac{E_3\sqrt{Z_2^2 - E_1^2}}{\sqrt{r(r-2)}}, \quad S^r= \frac{\sqrt{b(r)}}{rZ_2}(E_1 E_2 Z_1 + E_3 Z_2 Z_3) - \mathcal{E} E_1 - \frac{(r-3) Z_1 E_1}{r^2 Z_2}\sqrt{Z_2^2 - E_1^2}, \\
S^\theta = \frac{\sqrt{Z_2^2 - E_1^2}}{P_\varphi Z_2 r^\frac{5}{2}\sin \theta \sqrt{r-2}}\big( E_2 Z_1 (r-3)\sqrt{Z_2^2 - E_1^2} - Z_1(E_2^2 + E_3^2)\sqrt{r(r-2)} + r^2 \mathcal{E} E_2 Z_2 \big), \\
S^\varphi = \frac{E_3 P_\varphi \sqrt{Z_2^2 - E_1^2}}{r^{\frac{5}{2}} \sqrt{r -2} Z_2(E_2^2 + E_3^2)} \big( \mathcal{E} r^2 Z_2 + Z_1(r-3) \sqrt{Z_2^2 - E_1^2} \big),
\end{gathered}
$$
where in the last two components we have explicitly taken into account relations  \eqref{eqth} 
and the relation $P_r=0$, which is valid for relative equilibria. In addition, using the same relations, the linear momentum $p_\alpha$  can be represented as follows:
$$
\begin{gathered}
p_t = \frac{Z_1}{rZ_2}\sqrt{Z_2^2 - E_1^2} - \mathcal{E}, \quad p_r=0, \quad 
p_\theta = \frac{E_3}{ P_\varphi \sin \theta} \sqrt{b(r)(Z_2^2 - E_1^2)}, \\
p_\varphi = P_\varphi  - \frac{E_1^2}{P_\varphi} - \frac{E_2}{P_\varphi}\sqrt{b(r)(Z_2^2 - E_1^2)}.
\end{gathered}
$$
For the symmetric circular orbits, the linear momentum is parallel to the velocity ($p^\alpha \parallel u^\alpha$) and the following relations hold: 
$$
\mathcal{N}_0^\pm:  \  S^t=0, \ S^r=0, \ S^\theta \neq 0, \ S^t=0.
$$
As a consequence, for  $\mathcal{N}_0^\pm$ the total angular momentum is parallel to the angular momentum of the  test body.            
                
For the asymmetric circular orbits, it turns out that $p^\alpha \nparallel u^\alpha$ and the following relations hold:

$$
\begin{gathered}
\mathcal{N}_1^\pm: \ S^t\neq 0,  \ S^r=0, \ S^\theta \neq 0, \ S^\varphi \neq 0, \\
\mathcal{N}_2^\pm: \ S^t= 0, \ S^r \neq 0, \ S^\theta \neq 0, \ S^\varphi =0.
\end{gathered}
$$
From the fact that, in addition to $S^\theta$ there are other nonzero components of the Pauli-Lubanski vector it follows that, for the asymmetric circular orbits, the angular momentum of the test body is not parallel to the total angular momentum.

\subsection{Stability analysis}
\label{SectionStabilityAnalysis}

Let us analyze the stability of the asymmetric relative equilibria together with $\Sigma^+_0$.
Note that the pairs of relative equilibria are related to each other by the symmetry
transformation:
$$
 \mathcal{I}^{(r)}\big(\mathcal{N}_1^+\big)=\mathcal{N}_1^-, \quad \mathcal{I}^{(r)}\big(\mathcal{N}_2^+\big)=\mathcal{N}_2^-,
$$
and hence the equilibrium points from each pair are of the same type.

 The characteristic polynomial of linearization of the Hamiltonian system \eqref{eq_EZ}
 in a neighborhood of the equilibrium points on the invariant manifold \eqref{eq_Tred} has the form
$$
p(\lambda)=\lambda^2(\lambda^4 + A \lambda^2 + B),
$$
where two zero eigenvalues are related to the additional integrals \eqref{eq_Co_F}.
Depending on the character of the eigenvalues, four types of singular points are possible:
\begin{itemize}
\item[---] center-center ($\lambda_{1,2}=\pm i \alpha$, $\lambda_{3,4}=\pm i \beta$);
\item[---] saddle-center ($\lambda_{1,2}=\pm i \alpha$, $\lambda_{3,4}=\pm\beta$);
\item[---] saddle-saddle ($\lambda_{1,2}=\pm \alpha$, $\lambda_{3,4}=\pm \beta$);
\item[---] focus-focus ($\lambda_{1,2,3,4}=\pm \alpha \pm i \beta$),
\end{itemize}
where $\alpha$ and $\beta$ are real numbers. We note that it is only the singular point
of center-center type that is Lyapunov
stable.

For the equilibrium points $\mathcal{N}_1^\pm$ the coefficients of the characteristic polynomial $p(\lambda)$ have the form
 $$
 \begin{gathered}
A=\frac{(\ell -27 \mathcal{E}^2)^3}{81(\ell +54\mathcal{E}^2)^2} + \frac{4c^4 - 297c^2 - 1458}{27(\ell + 54 \mathcal{E}^2)^2}\mathcal{E}^2 -
\frac{3(4 c^2 - 27)}{(\ell + 54 \mathcal{E}^2)^2}(3\mathcal{E})^6 - \frac{2c^6 - 81 c^4 + 9^5\mathcal{E}^4c^2 -54^3}{(c^2 - 27)^2(\ell +54\mathcal{E}^2)^2}9\mathcal{E}^4, \\
B=\frac{\mathcal{E}^4(3 \mathcal{E}^2 - 1)}{(c^2 - 27)^2(\ell +54\mathcal{E}^2)^2}\big((c^2 - 27)^3 + 243 \mathcal{E}^2 c^2 (2c^2 + 243\mathcal{E}^2 - 135) \big).
\end{gathered}
 $$
For the equilibrium points $\mathcal{N}_2^\pm$ the coefficients of the characteristic 
polynomial $p(\lambda)$ are not given by explicit relations. For this reason, some of their 
numerical 
values are presented in Table \ref{table1}.

\begin{figure*}[!ht]
\begin{center}
\includegraphics[scale=0.68]{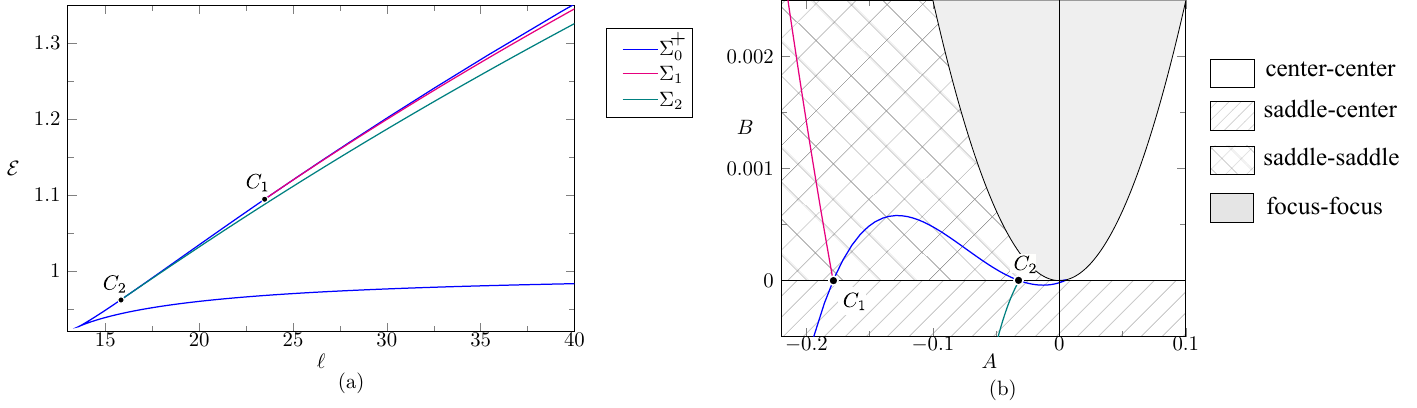}
\end{center}
\caption{ Curves $\Sigma_0^+$, $\Sigma_1$ and $\Sigma_2$ on the plane $(\ell, \mathcal{E})$
and the corresponding equilibrium points on the plane of the coefficients of the
characteristic polynomial for the fixed $c=0.65$, where $A$ and
$B$ are the coefficients of the characteristic polynomial defined in Section \ref{SectionStabilityAnalysis}.}
\label{fig8}
\end{figure*}

\begin{figure*}[!ht]
\begin{center}
\includegraphics[scale=0.68]{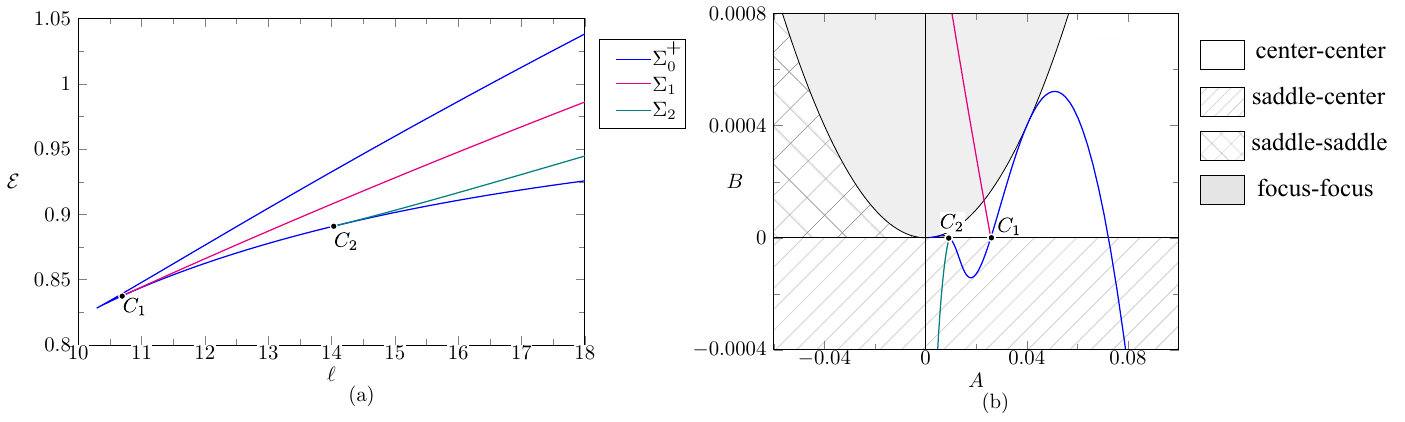}
\end{center}
\caption{Curves $\Sigma_0^+$, $\Sigma_1$ and $\Sigma_2$ on the plane $(\ell, \mathcal{E})$
and the corresponding equilibrium points on the plane of the coefficients of the
characteristic polynomial for the fixed $c=1.5$,  where $A$ and
$B$ are the coefficients of the characteristic polynomial defined in Section \ref{SectionStabilityAnalysis}.}
\label{fig9}
\end{figure*}

To analyze the stability of the equilibrium points, we will construct for the fixed
$c={\rm const}$ the curve $\Sigma_0^+$ together with $\Sigma_1$ and $\Sigma_2$ on the plane
$(\ell, \mathcal{E})$. In addition, we will show the corresponding equilibrium points
on the plane of the coefficients of the characteristic polynomial $(A,B)$ as curves.
This will allow us to define the type of equilibrium points (for details, see
\cite{BolsinovConleyIndex}) since regions with different types of singular
points on the plane of the coefficients are separated by the parabola $A^2 - 4B=0$ and the
straight line $B=0$ (see Figs. \ref{fig8}b and \ref{fig9}b). Both figures
provide a clear description of the types of equilibrium points and their
bifurcations.

We will distinguish three cases:
\begin{itemize}
\item [(i)] $c\in (0, c_1^*)$, where $c_1^*\approx 0.98878$;
\item [(ii)] $c\in (c_1^*, c_2^*)$, where $c_2^*\approx 1.3234$;
\item[(iii)] $c>c_2^*$.
\end{itemize}
 The typical form of the curves $\Sigma_0^+$, $\Sigma_1$ and $\Sigma_2$ for case (i)
 is shown in Fig. \ref{fig8}a. As we see, asymmetric circular orbits arise from $\Sigma_0^+$
 at points $C_1$ and $C_2$. The types of singular points for this case are shown in
 Fig.\ref{fig8}b, in which the type of the singular point $\Sigma_0^+$, except for the cusp,
 changes at points $C_1$ and $C_2$.

 The value $c=c_1^*$ corresponds to the case where point $C_2$ coincides with the cusp
 on the curve $\Sigma_0^+$, and then, as $c$ is increased, case (ii) arises, in which
 asymmetric circular orbits are also unstable. Then, for $c=c_2^*$ the other point $C_1$
 coincides with the cusp on the curve $\Sigma_0^+$. As $c$ is increased further,
 we obtain case (iii). For this case, the typical form of the curves $\Sigma_0^+$, $\Sigma_1$
 and $\Sigma_2$ and the types of singular points are shown in Fig. \ref{fig9}.
 As we see, in this case circular orbits $\mathcal{N}_1^\pm$ can be of center-center type
 and hence be stable.

The cusp on the curve $\Sigma_0^+$ corresponds to innermost stable circular orbits of spinning test particles \cite{Jefremov2015}. For this cusp, we denote the value of the radial coordinate $r$ by
$r_{ISQO}$.
Cases (ii) and (iii) imply that the angular momentum of the test body
has an extremely large value. As shown above, in these cases the symmetric circular orbits
become unstable with respect to asymmetric perturbations (i.e., perturbations that do not lie 
on ${\rm Fix } \ \mathcal{I}^{(r)}$) when $r>r_{ISQO}$. We describe these regions.
In case (ii), the symmetric circular orbits are unstable on the  segment of the curve 
$\Sigma_0^+$ between point $C_2$ and the cusp.
In case (iii), the symmetric circular orbits are unstable on the segment of the curve 
$\Sigma_0^+$ lying between points $C_1$ and $C_2$.

\begin{table}[h]
\centering
\caption{Numerical values of the points of the curves in Fig. \ref{fig8} and 
Fig. \ref{fig9}, which correspond to the equilibrium points $\mathcal{N}_2^\pm$,
where $A$ and
$B$ are the coefficients of the characteristic polynomial defined in Section \ref{SectionStabilityAnalysis}.}
\begin{tabular}{|c|c|c|c|c|}
\hline
$c$& $A$ & $B$ & $\ell$ & $\mathcal{E}$ \\
\hline
\ & -0.0341255 & -0.0000493 & 20.805948 & 1.044942 \\
\ & -0.0367417 & -0.0001110 & 25.015007 & 1.111712 \\
0.65 \ & -0.0395377 & -0.0001866 & 30.001444 & 1.186802\\
\ & -0.0425284 & -0.0002794 & 35.048713 & 1.258738 \\
\ & -0.0457302 & -0.0003931 & 39.916950 & 1.324690\\
\hline
\ & 0.0081308 & -0.0000548 & 14.086260 & 0.891579\\
\ & 0.0073306 & -0.0001081 & 15.006024 & 0.903201\\
1.5 & 0.0062743 & -0.0002031 & 16.001545 & 0.916641\\
\ & 0.0054216 & -0.0003049 & 17.004360 & 0.930605\\
\ & 0.0048092 & -0.0003962 & 17.991796 & 0.944558\\
\hline
\end{tabular}
\label{table1}
\end{table}

 \begin{rem}
Throughout this section, for the above-mentioned curves we have assumed $c \leqslant 2$.
If $c>2$, then the curves $\Sigma_0^\pm$, $\Sigma_1$ and $\Sigma_2$ can
bifurcate or disappear. For example, the radicand in $\ell^*$ can become negative.
\end{rem}

\section{Poincar\'{e} map}
\label{sectionMap}
Consider the trajectories of the reduced system \eqref{eq_EZ} that are different
from the equilibrium points found in
Sec. \ref{section_eqP}. For this we make use of a Poincar\'{e} map. Below we describe
the procedure for constructing it.
\begin{itemize}
\item[1.] As the secant of the vector field of system \eqref{eq_EZ} we choose a
plane that is given by the condition
$
E_1=0.
$
\item[2.] We will construct the Poincar\'{e} map in the variables $(r, P_r)$.
As for the other five variables, we will define them from the conditions defining
the invariant manifold \eqref{eq_Tred}, from the Casimir functions \eqref{eq_Co_F} and
from the integral corresponding to the mass \eqref{eq_mmm}.
The values on the plane $(r, P_r)$ for which there are no real roots will be
shown as gray points, which will be spaced equally apart.
\item[3.] For the initial conditions described above, we numerically integrate
system \eqref{eq_EZ}. At the same time, using the Henon method \cite{Henon1982}, we will
find intersections of the trajectories with the chosen section, on which
the quantity $\dfrac{dE_1}{d\tau}$ has the same sign.
Showing the values of the coordinates on the plane $(r, P_r)$, we obtain
a two-dimensional Poincar\'{e} point map.
\end{itemize}
\begin{figure}[!ht]
\begin{center}
\includegraphics[scale=0.9]{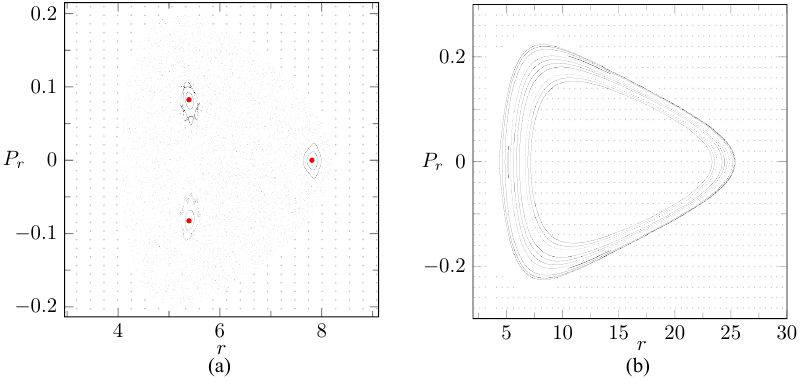}
\end{center}
\caption{A Poincar\'{e} map for fixed parameter values: (a) $c=1.4$, $\ell=16$, 
$E=0.92292941$, $m=1$ and $\mu=1$, which qualitatively reproduces the map from Ref. 
\cite[Fig. 4f]{Suzuki1997} and (b) $c=0.2$, $\ell=14.44$, $E=0.97$, $m=1$ and $\mu=1$, which qualitatively 
reproduces the map from Ref. \cite[Fig. 1]{Zelenka2020} }
\label{fig15}
\end{figure}
The procedure described above of constructing a Poincar\'{e} map differs from
that used in Ref. \cite{Suzuki1997}. The latter involves numerically solving part of the initial
MPD equations \eqref{eq_MP} in constructing a Poincar\'{e} map, whereas
the former involves numerically solving the reduced system \eqref{eq_EZ}, which has
two variables less. One variable has been eliminated using the Casimir function $C_\star=0$,
and the other, using the cyclic coordinate
of the symmetry field \eqref{eq_symu}.

The system was integrated using the explicit Dormand\,--\,Prince 8(7) method \cite{hairer1993solving} with automatic step-size control, for which the prescribed error tolerance was $10^{-12}$.
The integration was carried out under the condition that the conservation error of the 
system's additional integrals (the Hamiltonian, Casimir functions) does not exceed $10^{-9}$.
In this case, the Poincar\'{e} map for the parameter values from Refs. \cite{Suzuki1997} 
and \cite{Zelenka2020} is qualitatively reproduced (see Fig. \ref{fig15}) with the difference 
that, on the ordinate axis (the vertical axis), the authors of Refs. \cite{Suzuki1997} 
and  \cite{Zelenka2020} chose the velocity component of
the supporting curve $u^r$ and the linear momentum $p_r$, respectively. 
The authors of \ \cite{Suzuki1997} and \cite{Zelenka2020} chose the secant 
$\theta=\dfrac{\pi}{2}$, which, according to \eqref{eqth}, leads to the secant $E_1=0$ chosen 
in this paper. 

Nonetheless, in this case the Poincar\'{e} map for the parameter values from Ref. \cite{Suzuki1997}
is qualitatively reproduced (see Fig. \ref{fig15}) with the difference that, in
Ref. \cite{Suzuki1997}, the velocity component of the supporting curve $u^r$ (and not $P_r$)
and a different secant were chosen in the ordinate (vertical) axis.
\begin{figure*}[!ht]
\begin{center}
\includegraphics[scale=0.7]{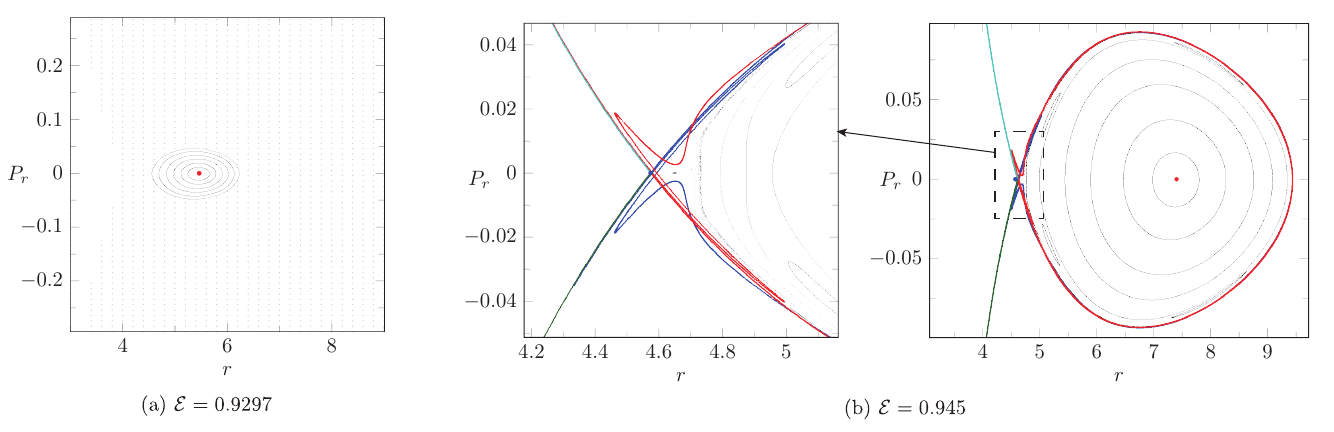}
\end{center}
\caption{A Poincar\'{e} map for the fixed $c=0.65$, $\ell=3.715$, $m=1$, $\mu=1$ and different
values of energy $\mathcal{E}$.}
\label{fig10}
\end{figure*}

We note that the trajectories corresponding to falling, $r\to 2\mu$, or those going to
infinity, $r\to + \infty$, can cease to cross the chosen secant. Therefore, a
Poincar\'{e} map can be used only in the case where
the recurrence of the trajectories is observed.
In Ref. \cite{Suzuki1997}, an analog of the effective potential has been constructed and
analyzed to search for such trajectories.

In this case, we make use of the results of Sec. \ref{section_eqP} and choose
values of the integrals $c$, $\ell$ and $\mathcal{E}$ a little above the part of the curve
$\Sigma^+_0$ in which the equilibrium point is of center-center type. A typical
Poincar\'{e} map for small values of $c$ is shown in Fig. \ref{fig10}a.
As we see, in this case the map has a stable fixed point (shown in red), which is
surrounded by invariant curves. Next, as we increase energy $\mathcal{E}$, we
cross the part of the curve $\Sigma^+_0$ that contains equilibrium points of saddle-center
type. An unstable fixed point arises in the map above this curve. This point is shown in
blue in Fig. \ref{fig10}b. In Fig. \ref{fig10}b,  stable and  unstable branches of the asymptotic manifolds to this
unstable fixed point are numerically constructed (they are shown in different colors).
One pair of asymptotic manifolds splits, and the other pair of asymptotic manifolds  tends to the event
horizon in forward or backward time. Previously, such a situation arose
for the geodesics of the Schwarzschild\,--\,Melvin metric \cite{Bizyaev2024}.
As energy $\mathcal{E}$ is increased further, the invariant curves in the map
break down. After some value of energy is reached, the trajectories begin
to fall onto the event horizon. We note that here and in what follows we assume $\mathcal{E}<1$,
for if this is not satisfied, trajectories going to infinity arise in the reduced system.

\begin{figure*}[!ht]
\begin{center}
\includegraphics[scale=0.8]{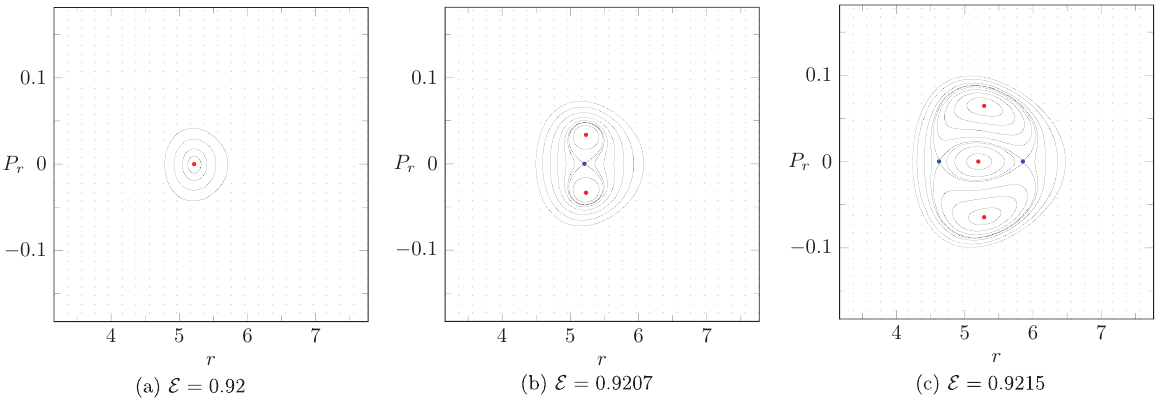}
\end{center}
\caption{A Poincar\'{e} map for the fixed $c=1.0$, $\ell=3.807$, $m=1$, $\mu=1$ and different
values of energy $\mathcal{E}$.}
\label{fig11}
\end{figure*}

\begin{figure*}[!ht]
\begin{center}
\includegraphics[scale=0.9]{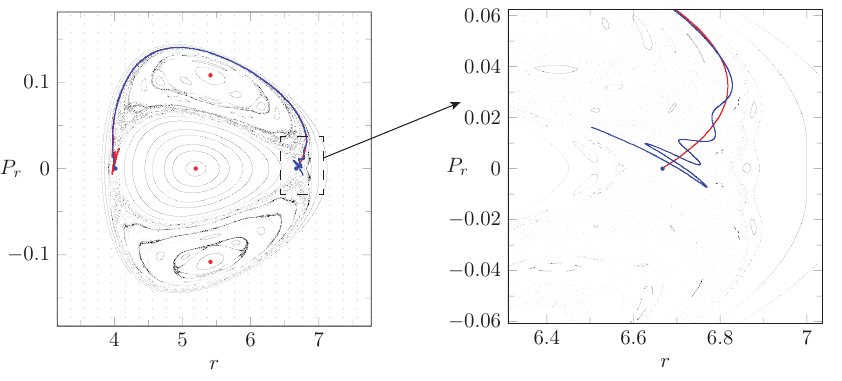}
\end{center}
\caption{A Poincar\'{e} map for the fixed $c=1.0$, $\ell=3.807$, $m=1$, $\mu=1$ and
$\mathcal{E}=0.9235$. }
\label{fig12}
\end{figure*}

\begin{figure*}[!ht]
\begin{center}
\includegraphics[scale=0.8]{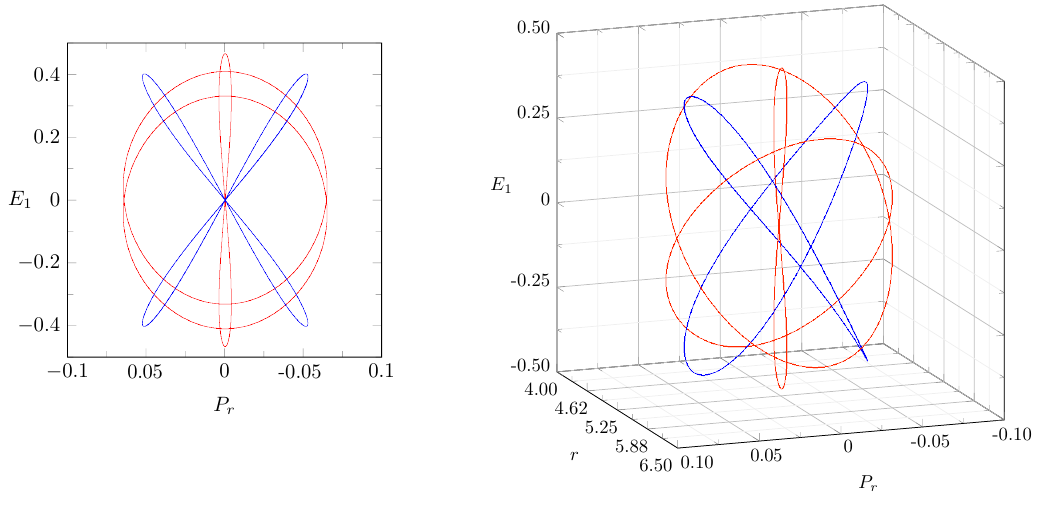}
\end{center}
\caption{Projection of the trajectories of the reduced system onto the plane $(P_r, E_1)$
and the space $(r, P_r, E_1)$ for the parameters and the initial conditions in a neighborhood
of fixed points
in the Poincar\'{e} map in Fig.\ref{fig11}c.}
\label{fig14}
\end{figure*}

For large values of $c$, a different scenario may take place. At first, the fixed point
in the map loses stability due to a supercritical pitchfork bifurcation
(see Fig. \ref{fig11}b), and then regains stability after a subcritical pitchfork
bifurcation (see Fig. \ref{fig11}c). Thus, five fixed points arise from one
fixed point. The corresponding trajectories of the reduced system are shown in
Fig. \ref{fig14}. Figure \ref{fig12} shows a Poincar\'{e} map and  asymptotic manifolds
constructed for different values of energy $\mathcal{E}$. As we see, asymptotic manifolds from different
fixed points in the map begin to intersect and heteroclinic trajectories arise.
This confirms the assumption made in Ref. \cite{Suzuki1997} that the observed chaos
is due to a ``heteroclinic tangle''. For these parameters, if we increase the energy further
after intersection with the curve $\Sigma_2$, the trajectories cease to cross the chosen secant
and begin to approach the event horizon.

\section{Conclusion}

This paper is concerned with the motion of a rotating test body in the Schwarzschild metric,
which is described by the MPD equations and the Tulczyjew condition.
This system admits additional integrals, making it possible to decrease the number of
equations in it, i.e., to perform a reduction.
This paper presents an explicit reduction after which the MPD equations
reduce to a Hamiltonian system with three degrees of freedom. This system has an invariant
manifold on which the required trajectories lie.
On this manifold, the problem reduces to a Hamiltonian system with two degrees of freedom.

Among the fixed points of this system there are the well-known \cite{Suzuki1997, Tod1976}
symmetric circular orbits. The total angular momentum of these trajectories is parallel to
the angular momentum of the test body. We find that the change of stability type of these
fixed points gives rise to asymmetric circular orbits. The total angular momentum of these
orbits is not parallel to the angular momentum of the body. If the test body is a black hole
or a neutron star, these solutions turn out to be unstable. We note that asymmetric
circular orbits have no analogs for geodesics since for them it holds that, as 
the spin magnitude
$c\to 0$, $\mathcal{E} \to +\infty$.

From the viewpoint of qualitative dynamical systems theory, the emergence of new equilibrium points can be explained as follows.
Indeed, the fixed point of the reduced system is simultaneously a critical point of
the Hamiltonian. As can be seen from Figs. \ref{fig8}b and \ref{fig9}b,
the symmetric circular orbits change their type if the
coefficient of the characteristic polynomial $B$ changes its sign. This is accompanied by
a change of the Morse index of the critical point. As a rule,
such a change of the index is accompanied by the emergence of new equilibrium points and helps
in finding them (for details, see, e.g.,
Ref. \cite{BolsinovConleyIndex}). Based on the results of Ref. \cite{Suzuki1997}, it appears that
such a situation also arises for symmetric circular orbits in a Kerr metric.
However, the analysis of the Kerr metric is made difficult by the fact
that, for this metric, the additional Carter integral \cite{Geoffrey2022} is not generalized
and the problem reduces to an analysis of
a Hamiltonian system with three degrees of freedom. Therefore, the search for asymmetric
circular orbits of a rotating test body for the Kerr metric is an interesting problem
in its own right.

In addition, this paper presents a Poincar\'{e} map constructed by numerically solving
the reduced system. Interestingly, this map exhibits two types of pitchfork bifurcations:
supercritical and subcritical. As a result, five fixed points arise in the Poincar\'{e} map.
Additionally, asymptotic manifolds are constructed for unstable fixed points which intersect
each other transversally. This confirms the assumption made in Ref. \cite{Suzuki1997}
that the observed chaos is due to a ``heteroclinic tangle''.

\bmhead{Acknowledgements}
The author extends appreciation to I.S. Mamaev for discussing the topics covered
in this paper and to the anonymous referees for their constructive criticisms, which 
have helped improve the original version of the manuscript.
This work was carried out within the framework of the
state assignment of the Ministry of Science and Higher
Education of Russia (FEWS-2024-0007).

\begin{appendices}
\section{Calculation of the Poisson bracket}\label{appendicesA}
We describe in more detail the calculation of the Poisson bracket between the reduced 
variables, which is presented in Section \ref{sectionRedQneq0}.

Denote the variables in the Poisson bracket \eqref{eq_J2} as follows: $\mathfrak{X} = (\vec{L}, \vec{M}, P_\alpha, y^\alpha)$. Then relations \eqref{eq_J2} define the following structure matrix:
$$
\begin{gathered}
\mathfrak{J}=
\left(
\begin{array}{ccc}
{\bf G}_{6\times 6} & {\bf O}_{6\times 4} & {\bf O}_{6\times 4} \\
 {\bf O}_{4\times 6} & {\bf O}_{4\times4} & -{\bf I}_{4 \times 4} \\
 {\bf O}_{4\times 6} & {\bf I}_{4 \times 4} & {\bf O}_{4\times4}
\end{array}
\right), \quad
{\bf G}_{6 \times 6} = \left(\begin{array}{cccccc}
0 & L_{3} & -L_{2} & 0 & M_{3} & -M_{2} \\

-L_{3} & 0 & L_{1} & -M_{3} & 0 & M_{1} \\

L_{2} & -L_{1} & 0 & M_{2} & -M_{1} & 0 \\

0 & M_{3} & -M_{2} & 0 & -L_{3} & L_{2} \\

-M_{3} & 0 & M_{1} & L_{3} & 0 & -L_{1} \\

M_{2} & -M_{1} & 0 & -L_{2} & L_{1} & 0
\end{array} \right),
\end{gathered}
$$
where ${\bf O}_{n \times m}$ is a zero matrix that has $n$ rows and $m$ columns and 
${\bf I}_{4 \times 4}$ is an identity matrix.

The structure matrix allows one to define the Poisson bracket (for details, 
see \cite{ArnoldKozlov, BorisovMamaev2018})
as follows:
\begin{equation}
\label{eq_FG}
\{ F, G\} = \mathfrak{J}^{ \hat{i} \hat{j}} \frac{\partial F}{ \partial \mathfrak{X}^{\hat{i}}} \frac{\partial G}{ \partial \mathfrak{X}^{\hat{j}}}.
\end{equation}
Here and below, the symbols $\hat{i}$ and $\hat{j}$ run over all phase variables 
$\mathfrak{X}$, i.e., from 1 to 14.

The equations of motion for the given Hamiltonian $H$ have the form
$$
\frac{d \mathfrak{X}^{\hat{i}} }{d \tau} = \{ \mathfrak{X}^{\hat{i}}, H \} = \mathfrak{J}^{ \hat{i} \hat{j}} \frac{\partial H}{ \partial \mathfrak{X}^{\hat{j}}}.
$$
Straightforward calculations show that these equations of motion have the form \eqref{eq_HamSys}.

Recall that, if $\{ F, G\} =0$, then the functions $F$ and $G$ are in involution. The 
Casimir function $C$ is in involution with all phase variables since it satisfies the equation
$$
 \mathfrak{J}^{ \hat{i} \hat{j}} \frac{\partial C}{ \partial \mathfrak{X}^{\hat{j}}} = 0.
$$
Let us pass to the calculation of the Poisson bracket between the reduced variables.
For this we calculate the Poisson bracket between the components of the vector
$$
\vec{E}= \left( L_1, \frac{L_3 P_\theta}{\sqrt{L_2^2 + L_3^2}} + L_2\frac{L_1 \cos \theta - P_\varphi}{\sin \theta \sqrt{L_2^2 + L_3^2}},
-\frac{L_2 P_\theta}{\sqrt{L_2^2 + L_3^2}} + L_3\frac{L_1 \cos \theta - P_\varphi}{\sin \theta \sqrt{L_2^2 + L_3^2}} \right).
$$
Using \eqref{eq_FG} and making direct calculations, we obtain
$$
\{ E_1, E_2 \} = E_3, \quad \{ E_1, E_3 \} = -E_2, \quad \{ E_2, E_3 \} = E_1,
$$
these commutation relations form the algebra $so(3)$
Next, we calculate the Poisson bracket between the components of the vector
$$
\vec{Z}=\left(-\frac{M_1 \sqrt{(\vec{L}, \vec{L})} }{\sqrt{L_2^2 + L_3^2}}, \sqrt{(\vec{L}, \vec{L})}, \frac{M_2L_3 - M_3L_2}{\sqrt{L_2^2 + L_3^2}} \right).
$$
Again using \eqref{eq_FG} and making direct calculations, we obtain
$$
\begin{gathered}
\{ Z_1, Z_2 \} = Z_3, \quad \{ Z_2, Z_3 \} = Z_1 +\frac{ (\vec{L}, \vec{M}) L_1}{\sqrt{(\vec{L}, \vec{L})(L_2^2 + L_3^2)}},\\
\{ Z_1, Z_3 \} = Z_2 + \frac{(\vec{L}, \vec{M})}{L_2^2 + L_3^2}\left( \frac{ \sqrt{(\vec{L}, \vec{L})}}{L_2^2 + L_3^2}(L_2M_2 + L_3M_3) - \frac{L_1 M_2}{ \sqrt{(\vec{L}, \vec{L})}} \right).
\end{gathered}
$$
On the zero level set of the Casimir function $C_\star= (\vec{L}, \vec{M})=0$, we obtain
commutation relations that form the algebra $so(2,1)$.
The variables $\vec{E}$ and $\vec{Z}$ do not explicitly depend on $r$, $P_r$ and hence are in involution 
with them. From the calculations made, we obtain the Poisson bracket \eqref{SP}.

\end{appendices}

\bibliography{Bizyaev}

\end{document}